\documentclass[10pt]{amsart}

\newtheorem{theorem}{Theorem}[section]
\newtheorem{proposition}[theorem]{Proposition}
\newtheorem{lemma}[theorem]{Lemma}
\newtheorem{corollary}[theorem]{Corollary}
\theoremstyle{definition}

\theoremstyle{remark} 
\numberwithin{equation}{section}

\begin{document}
\title[Multivariate generalized Linnik's probability densities]{Analytic and asymptotic properties of multivariate  generalized Linnik's probability densities}

\author{S. C. Lim}
\address{Faculty of Engineering, Multimedia University, Jalan Multimedia, 63100, Cyberjaya, Selangor.}
\email{sclim@mmu.edu.my}

\author{L. P. Teo}
\address{Faculty of Information Technology, Multimedia University, Jalan Multimedia, 63100, Cyberjaya, Selangor.}
\email{lpteo@mmu.edu.my}
\begin{abstract}This paper studies the   properties of the probability density function $p_{\alpha,\nu, n}(\mathbf{x})$ of the $n$-variate generalized Linnik distribution whose characteristic function $\varphi_{\alpha,\nu,n}(\boldsymbol{t})$ is given by
\begin{equation*}
\varphi_{\alpha,\nu,n}(\boldsymbol{t})=\frac{1}{\left(1+ \Vert\boldsymbol{t}\Vert^{\alpha}\right)^{\nu}}, \;\;\;\;\alpha\in (0,2], \;\nu>0,\;\boldsymbol{t}\in \mathbb{R}^n,
\end{equation*}where $\Vert\boldsymbol{t}\Vert$ is the Euclidean norm of $\boldsymbol{t}\in\mathbb{R}^n$. Integral representations of $p_{\alpha,\nu, n}(\mathbf{x})$ are obtained and used to derive the asymptotic expansions of $p_{\alpha,\nu, n}(\mathbf{x})$ when $\Vert\mathbf{x}\Vert\rightarrow 0$ and $\Vert\mathbf{x}\Vert\rightarrow \infty$ respectively. It is shown that under certain conditions which are arithmetic in nature,   $p_{\alpha,\nu, n}(\mathbf{x})$ can be represented   in terms of entire functions.

\end{abstract}
\maketitle

\section{Introduction}
In 1953, Linnik \cite{1} showed that the function
\begin{equation}\label{eq10_28_1}
\varphi_{\alpha}(t)=\frac{1}{1+|t|^{\alpha}}, \hspace{1cm}
\alpha\in (0,2], \;\;t\in \mathbb{R},
\end{equation}is the characteristic function of a symmetric probability density $p_{\alpha}(x)$. Since then, $p_{\alpha}(x)$ is known as Linnik's probability density. When $\alpha=2$, $p_2(x) $ is   the  probability density of the  Laplace  distribution. Therefore, $p_{\alpha}(x)$ is sometimes referred to as $\alpha$-Laplace probability density.
In \cite{11}, Klebanov, Maniya and Melamed introduced the concept of geometric strict stability and showed that the characteristic functions of the probability distributions that are geometrically strictly stable are up to a scaling factor given by:
\begin{equation}\label{eq10_28_2}
\varphi_{\alpha,\theta}(t)=\frac{1}{1+e^{-i\theta \;\text{sgn} \,t} |t|^{\alpha}}, \hspace{1cm}\alpha\in (0,2), \;\;|\theta|\leq \min\left(\frac{\pi\alpha}{2}, \pi-\frac{\pi\alpha}{2}\right).
\end{equation}The corresponding probability distribution is known as asymmetric Linnik distribution.  In \cite{8}, Pakes showed that the probability distributions known as generalized Linnik distributions which have characteristic functions
\begin{equation}\label{eq10_28_3}
\varphi_{\alpha, \nu, \theta}(t)=\frac{1}{\left(1+e^{-i\theta \;\text{sgn}\, t} |t|^{\alpha}\right)^{\nu}}, \hspace{1cm}\alpha\in (0,2), \;\;|\theta|\leq \min\left(\frac{\pi\alpha}{2}, \pi-\frac{\pi\alpha}{2}\right), \;\;\nu>0,
\end{equation}play an important role in some characterization problems of mathematical statistics.
The class of probability distributions with density functions $p_{\alpha, \nu, \theta}(x)$ have found some interesting properties and applications \cite{3, 4, 36, 5, 7, 6, 43,  13, 31, 38, 37, 14, 25, 28, 34, 39, 12, 16, 26, 27, 29, 24, 32, 33, 8, 9, 22,  35}.
In particular, they are good candidates to model financial data which exhibits high kurtosis and heavy tails.  In \cite{30} and \cite{2}, the probability distribution $p_{\alpha}(x)$ was independently generalized to   multivariate distributions. Anderson \cite{2} showed that the function
\begin{equation}\label{eq10_28_4}
\varphi_{\alpha, n, \Sigma}(\boldsymbol{t}) = \frac{1}{1+ \left[ \boldsymbol{t}' \Sigma \boldsymbol{t}\right]^{\frac{\alpha}{2}}}, \hspace{1cm}\alpha\in (0, 2], \;\; \boldsymbol{t}\in\mathbb{R}^n,
\end{equation}where $\Sigma$ is a positive definite $n\times n$ matrix, is a characteristic function of a $n$-variate probability distribution which he called the $n$-variate Linnik distribution. Properties and applications of this distribution have been studied in \cite{21, 42, 69, 40, 23, 12, 41, 16, 32}. Similar to the univariate case, the multivariate Linnik distribution is a subclass of multivariate geometric stable distributions \cite{32, 40}.

As the case of stable distributions, in general the probability density functions of the univariate and multivariate Linnik distributions do not have closed forms. In \cite{10}, Kotz, Ostrovskii and Hayfavi studied the analytic and asymptotic behaviors of the probability density function of the univariate symmetric Linnik distribution whose characteristic function is given by \eqref{eq10_28_1}. This work was generalized to symmetric homogeneous multivariate Linnik distribution with characteristic function   \eqref{eq10_28_4} with $\Sigma=\text{id}$ by Ostrovskii \cite{18}; and to asymmetric univariate Linnik distribution with characteristic function   \eqref{eq10_28_2} by Erdogan \cite{17}, and finally to  asymmetric univariate Linnik distribution with characteristic function \eqref{eq10_28_3} by Erdogan and Ostrovskii \cite{20, 19}. In the present work, we consider the generalization of the multivariate Linnik distribution whose characteristic function is given by
\begin{equation}\label{eq10_28_5}
\varphi_{\alpha,\nu, n, \Sigma}(\boldsymbol{t}) = \frac{1}{\left(1+ \left[ \boldsymbol{t}' \Sigma \boldsymbol{t}\right]^{\frac{\alpha}{2}}\right)^{\nu}}, \hspace{1cm}\alpha\in (0, 2], \;\nu>0,\; \boldsymbol{t}\in\mathbb{R}^n.
\end{equation}One can expect that this multivariate generalized Linnik distribution also plays an important role in some characterization problems in multivariate statistics as in the univariate case \cite{8}. The main goal of this paper is to investigate  the asymptotic and analytic properties of the probability density function of the multivariate generalized Linnik distribution. It suffices to restrict to the case where $\Sigma=\text{id}$, i.e., we only consider the class of multivariate symmetric generalized Linnik distribution with characteristic function given by
\begin{equation}\label{eq10_28_6}
\varphi_{\alpha,\nu, n}(\boldsymbol{t}) = \frac{1}{\left(1+ \Vert \boldsymbol{t}\Vert^{ \alpha }\right)^{\nu}}, \hspace{1cm}\alpha\in (0, 2], \;\nu>0,\; \boldsymbol{t}\in\mathbb{R}^n,
\end{equation}where $\Vert \boldsymbol{t}\Vert =\sqrt{\sum_{i=1}^n t_i^2}$ is the Euclidean norm of $\boldsymbol{t}$. We denote by $p_{\alpha,\nu, n}(\mathbf{x})$ the corresponding probability density function.

There are also other motivations to the problem we want to study here. Recall that the dual distribution $\hat{p}$ of a probability distribution $p$ is a probability distribution whose characteristic function (resp. probability density function) is up to a constant, the probability density function (resp. characteristic function) of $p$ \cite{15}. When $\alpha\in (0, 2]$ and $\nu>n /\alpha$, the dual distribution of the generalized Linnik distribution is called the   generalized Cauchy distribution \cite{44} whose probability density function is up to a constant, given by \eqref{eq10_28_6}. Therefore, $p_{\alpha,\nu, n}(\mathbf{x})/p_{\alpha, \nu, n}(\mathbf{0})$ is   the characteristic function of the generalized Cauchy distribution.  For the special case where $\alpha=2$ and $\nu>n/2$, the generalized Cauchy distribution is the well-known  Student's $t$ distribution  which has a lots of applications (see e.g. \cite{45} and references therein).

From the perspective of stochastic processes, Gneiting and Schlater \cite{46} introduced a new class of stationary stochastic processes called Gaussian field with generalized Cauchy covariance whose covariance function is given by \eqref{eq10_28_6}. This stochastic model has two parameters $\alpha$ and $\nu$ which can give separate characterizations of the fractal dimension and Hurst effect. It has found applications in many modeling problems \cite{47, 53, 48, 52, 49,   55, 50, 51}. In this context, the function under investigation $p_{\alpha,\nu, n}(\mathbf{x})$ appears as the spectral density function of the stochastic model \cite{54}. On the other hand, when $\alpha=2$ and $\nu> n/2$, up to a constant, \eqref{eq10_28_6} appears as the spectral density of the Whittle-Mat$\acute{\text{e}}$rn random field \cite{58, 56, 57}. This class of random fields has found wide applications especially in geostatistics \cite{15, 62, 63, 64, 59,   61, 65, 66, 60}. In \cite{67}, we generalized the Whittle-Mat$\acute{\text{e}}$rn random field to a random field whose spectral density is given by \eqref{eq10_28_6} with $\alpha\in (0, 2]$, $\nu >n/\alpha$ and showed that it can provide a more flexible model for wind speed data. The function $p_{\alpha,\nu, n}(\mathbf{x})$ we want to study in this paper then becomes the covariance function of the generalized Whittle-Mat$\acute{\text{e}}$rn random field.

From the brief discussion above, one notes that the function $p_{\alpha,\nu,n}(\mathbf{x})$ appears in various places and plays different roles, stretching from the probability density function of generalized Linnik distribution and spectral density function of Gaussian field with generalized Cauchy covariance to the characteristic function of the generalized Cauchy distribution and the covariance function of generalized Whittle-Mat$\acute{\text{e}}$rn field. Thus, a detailed study of the analytic and asymptotic properties of $p_{\alpha,\nu,n}(\mathbf{x})$ is of strong interest and significant importance, particularly in view of its potential applications in physics, internet traffic and financial time series analysis and modeling. This is the task undertaken in the present work.

\section{Multivariate generalized Linnik distribution}   Gneiting and Schlather \cite{46} asserted that using similar  arguments as  of \cite{68} and references therein, one can show that the function $\varphi_{\alpha, \nu, n}(\boldsymbol{t})$  \eqref{eq10_28_6} is a covariance function of a stationary Gaussian random field if and only if $\alpha\in (0, 2]$ and $\nu>0$. This implies that when $\alpha\in (0, 2]$ and $\nu>0$, $p_{\alpha,\nu, n}(\mathbf{x})\geq 0$ and therefore is also the probability density function of a distribution. In this section, we provide another argument to show that $\varphi_{\alpha, \nu, n}(\boldsymbol{t})$  \eqref{eq10_28_6} is the characteristic function of a probability distribution. The following is a generalization of the result of Devroye \cite{7} to multivariate case.

\begin{proposition}\label{Thm1}
Given $\alpha\in (0, 2]$, $\nu>0$ and $\Sigma$ a positive definite matrix of rank $n$, let $\mathbf{S}_{\alpha,  n, \Sigma}$ be a symmetric multivariate stable  random variable with characteristic function $\exp\left(-\left[ \boldsymbol{t}' \Sigma \boldsymbol{t}\right]^{\frac{\alpha}{2}}\right)$ and let $U_{\nu}$ be an independent univariate gamma random variable with probability density function $$\frac{u^{\nu-1}}{\Gamma(\nu)}e^{-u}, \hspace{1cm} u\geq 0.$$Then the characteristic function of the random vector $\mathbf{X}_{\alpha,\nu,n, \Sigma}= U_{\nu}^{\frac{1}{\alpha}}\mathbf{S}_{\alpha, n, \Sigma}$ is given by \eqref{eq10_28_5}.
\end{proposition}
\begin{proof}
\begin{equation*}\begin{split}
E\left( e^{i\boldsymbol{t}.\mathbf{X}}\right)=&E\left(E\left( \left.e^{iU_{\nu}^{\frac{1}{\alpha}}\boldsymbol{t}.\mathbf{S}_{\alpha,n, \Sigma}}\right| U_{\nu}\right)\right)\\
=&E\left( \exp\left(-U_{\nu}\left[ \boldsymbol{t}' \Sigma \boldsymbol{t}\right]^{\frac{\alpha}{2}}\right)\right)\\
=&\frac{1}{\Gamma(\nu)}\int_0^{\infty} u^{\nu-1} \exp\left\{-u\left(1+\left[ \boldsymbol{t}' \Sigma \boldsymbol{t}\right]^{\frac{\alpha}{2}}\right)\right\}du\\
=&\frac{1}{\left(1+\left[ \boldsymbol{t}' \Sigma \boldsymbol{t}\right]^{\frac{\alpha}{2}}\right)^{\nu}}.
\end{split}\end{equation*}
\end{proof}

We call the random vector $\mathbf{X}_{\alpha,\nu,n}$ whose characteristic function is given by \eqref{eq10_28_6} an $(\alpha, \nu, n)$ Linnik random vector.
As in the introduction, denote by $p_{\alpha,\nu,n}(\mathbf{x})$ its probability density function. In other words, $p_{\alpha,\nu,n}(\mathbf{x})$ is the unique function such that
\begin{equation}\label{eq11_12_7}
\varphi_{\alpha, \nu, n}(\boldsymbol{t})=\frac{1}{\left(1+\Vert \boldsymbol{t}\Vert^{ \alpha }\right)^{\nu}}=\int_{\mathbb{R}^n}e^{i\mathbf{x}. \boldsymbol{t}}p_{\alpha,\nu,n}(\mathbf{x})d^n\mathbf{x}.
\end{equation}Since $\varphi_{\alpha, \nu, n}(\boldsymbol{t})$ is a radial function,  $p_{\alpha,\nu,n}(\mathbf{x})$ is also a radial function. Denoting by $q_{\alpha, \nu, n}(r)$ the function $$q_{\alpha,\nu, n}(r)= \left.p_{\alpha,\nu,n}( \mathbf{x} )\right|_{\Vert \mathbf{x}\Vert=r}.$$
Then Schoenberg's formula gives
\begin{equation}\label{eq10_29_4}
\frac{1}{\left(1+\Vert \boldsymbol{t}\Vert^{ \alpha }\right)^{\nu}}=(2\pi)^{\frac{n}{2}}\Vert t\Vert^{\frac{2-n}{2}}\int_0^{\infty} J_{\frac{n-2}{2}}(r\Vert t\Vert)q_{\alpha,\nu, n}(r)r^{\frac{n}{2}}dr,
\end{equation}
where $J_{\mu}(z)$ is the   Bessel function of the first kind. In this paper, we want to study the properties of the function $q_{\alpha,\nu, n}(r)$.

Proposition \ref{Thm1} gives a relation between the probability density functions of the generalized Linnik distribution and the symmetric stable distribution. More precisely, let $S_{\alpha, n}(\mathbf{x})$ be the probability density function of the symmetric stable distribution with characteristic function $\exp(-\Vert \boldsymbol{t}\Vert^{\alpha})$, i.e.
\begin{equation}\label{eq11_12_5}
e^{-\Vert \boldsymbol{t}\Vert^{\alpha}}=\int_{\mathbb{R}^n} e^{i\boldsymbol{t}.\mathbf{x}}s_{\alpha, n}(\mathbf{x})d^n\mathbf{x},
\end{equation}
or equivalently,
\begin{equation}\label{eq11_12_6}
s_{\alpha, n}(\mathbf{x})=\frac{1}{(2\pi)^n}\int_{\mathbb{R}^n} e^{-i\boldsymbol{t}.\mathbf{x}}e^{-\Vert \boldsymbol{t}\Vert^{\alpha}}d^n \boldsymbol{t}=\frac{\Vert\mathbf{x}\Vert^{\frac{2-n}{2}}}{(2\pi)^{\frac{n}{2}}}\int_0^{\infty}J_{\frac{n-2}{2}}(t\Vert\mathbf{x}\Vert)e^{-t^{\alpha} }t^{\frac{n}{2}}dt.
\end{equation}Then
\begin{proposition}\label{pro13}
Let $\alpha\in (0, 2]$ and $\nu>0$. The probability density function $p_{\alpha,\nu, n}(\mathbf{x})$ has a representation
\begin{equation}\label{eq10_29_5}
p_{\alpha,\nu, n}(\mathbf{x})=q_{\alpha,\nu, n}(\Vert\mathbf{x}\Vert)=\frac{\alpha}{\Gamma(\nu)}\int_0^{\infty} u^{\alpha\nu-n-1}e^{-u^{\alpha}}\tilde{s}_{\alpha, n}\left(\frac{\Vert\mathbf{x}\Vert}{u}\right)du,
\end{equation}where $\tilde{s}_{\alpha, n}(r) =\left. s_{\alpha, n}(\mathbf{x})\right|_{\Vert \mathbf{x}\Vert=r}$.
\end{proposition}The properties of the function $\tilde{s}_{\alpha, n}(r)$ have been studied in \cite{72}. In particular, it was shown that
\begin{equation}\begin{split}\label{eq11_12_2}
\tilde{s}_{\alpha, n}(r)=&O( r^0) \hspace{2cm} \text{as}\;\; r\rightarrow 0,\\
\tilde{s}_{\alpha, n}(r)=&O\left( r^{-n-\alpha}\right) \hspace{1cm} \text{as}\;\; r\rightarrow \infty.
\end{split}\end{equation}These show that \eqref{eq10_29_5} is indeed well-defined for all $\alpha\in (0, 2]$, $\nu>0$ and $r>0$.
\begin{proof}
Making a change of variable, we can rewrite the right hand side of \eqref{eq10_29_5} as
\begin{equation}\label{eq11_12_1}
\frac{1}{\Gamma(\nu)}\int_0^{\infty} u^{\nu-\frac{n}{\alpha}-1}e^{-u }\tilde{s}_{\alpha, n}\left(\frac{\Vert\mathbf{x}\Vert}{u^{\frac{1}{\alpha}}}\right)du.
\end{equation}Taking the Fourier transform of \eqref{eq11_12_1} and using the definition  of $\tilde{s}_{\alpha, n}(r)$ \eqref{eq11_12_5}, we obtain
\begin{equation*}\begin{split}
&\frac{1}{\Gamma(\nu)}\int_{\mathbb{R}^n}e^{i\boldsymbol{t}. \mathbf{x}}\int_0^{\infty} u^{\nu-\frac{n}{\alpha}-1}e^{-u }\tilde{s}_{\alpha, n}\left(\frac{\Vert\mathbf{x}\Vert}{u^{\frac{1}{\alpha}}}\right)du d^n\mathbf{x}\\
=&\frac{1}{\Gamma(\nu)}\int_0^{\infty}\int_{\mathbb{R}^n}e^{iu^{\frac{1}{\alpha}}\boldsymbol{t} . \mathbf{x}}u^{\nu-1}e^{-u}\tilde{s}_{\alpha, n}\left( \Vert\mathbf{x}\Vert \right)  d^n\mathbf{x}du\\
=&\frac{1}{\Gamma(\nu)}\int_0^{\infty} u^{\nu-1}e^{-u\left(1+\Vert \boldsymbol{t}\Vert^{\alpha}\right)} du\\
=&\frac{1}{\left(1+\Vert \boldsymbol{t}\Vert^{\alpha}\right)^{\nu}}.
\end{split}\end{equation*}This proves the assertion of Proposition \ref{pro13}. The interchange of the order of integrations can be justified using \eqref{eq11_12_2}.

\end{proof}As is observed in \cite{72}, the identity (see e.g. \cite{71})
\begin{equation*}
\frac{d}{dz}J_{\mu}(z)=\frac{\mu}{z}J_{\mu}(z)-J_{\mu+1}(z)
\end{equation*}and \eqref{eq11_12_6} imply that
\begin{equation*}
\tilde{s}_{\alpha, n}'(r) =-2\pi r \tilde{s}_{\alpha, n+2}(r).
\end{equation*}It follows from \eqref{eq10_29_5} that
\begin{equation}\label{eq10_29_6}
q_{\alpha, \nu, n}'(r)=-2\pi r q_{\alpha, \nu, n+2}(r).
\end{equation}
Since $q_{\alpha, \nu, n+2}(r)\geq 0$, we deduce from \eqref{eq10_29_6} that

\begin{proposition}\label{pro14}
$q_{2,\nu, n}(r)$, $r\in \mathbb{R}_+$, is a monotonically decreasing function of $r$. \end{proposition}
Now we want to find other representations for $q_{\alpha,\nu,n}(r)$. A naive application of Fourier inversion formula to \eqref{eq11_12_7} gives
\begin{equation}\label{eq10_29_1}\begin{split}
p_{\alpha,\nu, n}(\mathbf{x}) =& \frac{1}{(2\pi)^n}\int_{\mathbb{R}^n} \frac{e^{-i\mathbf{x}. \boldsymbol{t}}}{\left(1+ \Vert \boldsymbol{t}\Vert^{ \alpha }\right)^{\nu}}d^n\boldsymbol{t}\\
=&\frac{\Vert \mathbf{x}\Vert^{\frac{2-n}{2}}}{(2\pi)^{\frac{n}{2}}}\int_0^{\infty} \frac{J_{\frac{n-2}{2}}(\Vert \mathbf{x}\Vert t)}{(1+ t^{\alpha})^{\nu}}t^{\frac{n}{2}}dt.\end{split}
\end{equation}However, since (see e.g. \cite{71})
\begin{equation*}
J_{\mu}(z)\sim \sqrt{\frac{2}{\pi z}}\cos\left( z-\frac{\pi \mu}{2}-\frac{\pi}{4}\right) \hspace{0.5cm}\text{as}\;\;z\rightarrow \infty,
\end{equation*}\eqref{eq10_29_1} is only valid when $\alpha\nu>\frac{n-1}{2}$. In \cite{54}, we proved the following:

\begin{proposition}\label{Pro2}
For $\alpha\in (0, 2)$ and $\nu>0$,  $q_{\alpha,\nu, n}(r) $ has the following integral representation
\begin{equation}\label{eq10_29_2}\begin{split}
q_{\alpha,\nu, n}(r) =&-\frac{r^{\frac{2-n}{2}}}{2^{\frac{n-2}{2}}\pi^{\frac{n+2}{2}}}\text{Im}\,\int_0^{\infty}
\frac{K_{\frac{n-2}{2}}(r
u)}{\left(1+e^{\frac{i\pi\alpha}{2}}u^{\alpha}\right)^{\nu}}u^{\frac{n}{2}}du.\end{split}
\end{equation}For $\alpha=2$ and $\nu>0$, $q_{2,\nu, n}(r) $ is given explicitly by
\begin{equation}\label{eq10_29_3}
q_{2,\nu, n}(r)= \frac{r^{\nu-\frac{n}{2}}}{2^{\frac{n}{2}+\nu-1}\pi^{\frac{n}{2}}\Gamma(\nu)}K_{\frac{n}{2}-\nu}
(r).
\end{equation}Here $K_{\mu}(z)$ is the modified Bessel function of the second kind.\end{proposition}
Since $K_{\mu}(z)$ is an infinitely differentiable function and \begin{equation}\label{eq11_04_2}\begin{split}K_{\mu}(z)=& O\left( |z|^{-|\mu|}\log\frac{1}{|z|}\right) \hspace{1.5cm}\text{as}\;\; z\rightarrow 0,\\
K_{\mu}(z) =  &O(e^{-z})\hspace{1cm}\text{as}\;\; z\rightarrow \infty,\end{split}\end{equation} (see e.g. \cite{71}), Proposition \ref{Pro2} shows that $q_{\alpha,\nu, n}(r)$ is an infinitely differentiable function of $r$ when $r\in (0, \infty)$.

When $\nu=1$, \eqref{eq10_29_2} has been proved by Ostrovskii \cite{18}. When $n=1$, it was proved in \cite{1, 10, 20, 19}. Both \eqref{eq10_29_2} and \eqref{eq10_29_3} can be proved by showing that they satisfy  \eqref{eq10_29_4}. However, to motivate these formulas, one can start from \eqref{eq10_29_1} when $\alpha\nu> \frac{n-1}{2}$. In this case,    the formula \eqref{eq10_29_1} can be evaluated explicitly when $\alpha=2$ which gives \eqref{eq10_29_3}. When $\alpha\in (0,2)$, we use the fact that  (see e.g. \cite{71}) $$H_{\mu}^{(1)}(z)=J_{\mu}(z)+iN_{\mu}(z),$$and \begin{equation*}
K_{\mu}(z)=\frac{i\pi}{2}e^{\frac{i\mu\pi}{2}}H_{\mu}^{(1)}(iz),
\end{equation*} where $H_{\mu}^{(1)}(z)$ is the Hankel's function of the first kind and
$N_{\mu}(z)$ is the Bessel function of the second kind or
called the Neumann function. A change of contour of integration  from the positive real axis to the positive imaginary axis gives \eqref{eq10_29_2}. For details, see \cite{54}.

Using Proposition \ref{Pro2}, we can prove a result stronger than the result of Proposition \ref{pro14}. When $\alpha=2$,  using the formula (see e.g. \cite{71})
 \begin{equation}\label{eq11_12_8}
 K_{\mu}(z)=\frac{\sqrt{\pi}\left(\frac{z}{2}\right)^{\mu}}{\Gamma\left(\mu+\frac{1}{2}\right)}\int_1^{\infty} e^{-yz}(y^2-1)^{\mu-\frac{1}{2}}dy,\hspace{0.5cm} \mu>-\frac{1}{2},\;z>0,
 \end{equation}we can rewrite \eqref{eq10_29_3} as
\begin{equation*}
q_{2,\nu, n} (r)= \frac{1}{2^{n-1}\pi^{\frac{n-1}{2}}\Gamma(\nu)\Gamma\left(\frac{n+1}{2}-\nu\right)}\int_1^{\infty} (y^2-1)^{\frac{n-1}{2}-\nu}e^{-ry}dy
\end{equation*}when $\nu <\frac{n+1}{2}$. When $\nu =\frac{n+1}{2}$, we have the explicit formula
\begin{equation*}
q_{2,\frac{n+1}{2}, n} (r)=\frac{e^{-r}}{2^{n}\pi^{\frac{n-1}{2}} \Gamma\left(\frac{n+1}{2} \right)}.
\end{equation*}
Therefore $(-1)^k q_{2,\nu, n}^{(k)} (r)\geq 0$ for all $k=0, 1, 2, \ldots$ and $q_{2,\nu, n} (r)$ is a completely monotonic function of $r\in \mathbb{R}_+$ when $\nu \leq \frac{n+1}{2}$. When $\alpha\in (0, 2)$, we rewrite \eqref{eq10_29_2} as
\begin{equation}\label{eq10_29_7} \begin{split}
q_{\alpha,\nu, n}(r) =&\frac{r^{\frac{2-n}{2}}}{2^{\frac{n-2}{2}}\pi^{\frac{n+2}{2}}}\int_0^{\infty}
\frac{\sin \left(\nu \text{arg}\left(1+e^{\frac{i\pi\alpha}{2}}u^{\alpha}\right)\right) }{\left\vert 1+e^{\frac{i\pi\alpha}{2}}u^{\alpha}\right\vert^{\nu}}K_{\frac{n-2}{2}}(r
u)u^{\frac{n}{2}}du.\end{split}
\end{equation}If $n=1$, this gives
\begin{equation}\label{eq10_29_8}
q_{\alpha,\nu, 1}(r) =\frac{1}{\pi}\int_0^{\infty}\frac{\sin \left(\nu \text{arg}\left(1+e^{\frac{i\pi\alpha}{2}}u^{\alpha}\right)\right) }{\left\vert 1+e^{\frac{i\pi\alpha}{2}}u^{\alpha}\right\vert^{\nu}}e^{-ru}du.
\end{equation}When $n\geq 2$, we can use \eqref{eq11_12_8} to transform \eqref{eq10_29_7} to
\begin{equation}\label{eq10_29_9}\begin{split}
q_{\alpha,\nu, n}(r) =&\frac{1}{2^{n-2}\pi^{\frac{n+1}{2}}\Gamma\left(\frac{n-1}{2}\right)}\int_0^{\infty}\frac{\sin \left(\nu \text{arg}\left(1+e^{\frac{i\pi\alpha}{2}}u^{\alpha}\right)\right) }{\left\vert 1+e^{\frac{i\pi\alpha}{2}}u^{\alpha}\right\vert^{\nu}}u^{n-1}\int_1^{\infty}(y^2-1)^{\frac{n-3}{2}}e^{-yru}dy du\\
=&\frac{1}{2^{n-2}\pi^{\frac{n+1}{2}}\Gamma\left(\frac{n-1}{2}\right)}\int_0^{\infty}e^{-yr}\left\{\int_0^y \frac{\sin \left(\nu \text{arg}\left(1+e^{\frac{i\pi\alpha}{2}}u^{\alpha}\right)\right) }{\left\vert 1+e^{\frac{i\pi\alpha}{2}}u^{\alpha}\right\vert^{\nu}}(y^2-u^2)^{\frac{n-3}{2}}udu\right\}dy.
\end{split}\end{equation}Since $\sin \left(\nu \text{arg}\left(1+e^{\frac{i\pi\alpha}{2}}u^{\alpha}\right)\right)\geq 0$ if $\alpha\nu\leq 2$, \eqref{eq10_29_8} and \eqref{eq10_29_9} imply that $q_{\alpha, \nu, n}(r)$ is completely monotonic if $\alpha\in (0,2)$ and $\alpha\nu \in (0, 2]$. The result is summarized in the following proposition.

\begin{proposition}
The function $q_{\alpha, \nu, n}(r)$ is completely monotonic if $\alpha\in (0,2)$ and $\alpha\nu\in (0,2]$ or $\alpha=2$ and $\nu\in \left(0, \frac{n+1}{2}\right]$.
\end{proposition}We remark that when $\alpha\in (0, 2)$ and $\nu=1$, this result has been proved in \cite{18}. When $n=1$, it has been shown in \cite{20, 19}.

\section{Asymptotic expansions of $q_{\alpha, \nu, n}(r)$}In this section, we derive the asymptotic expansions of the function $q_{\alpha, \nu, n}(r)$ when $r\rightarrow 0$ and $r\rightarrow \infty$ respectively. We discuss the cases $\alpha=2$ and $\alpha\in (0, 2)$ separately.

\subsection{$\alpha=2$} When $\alpha=2$, the asymptotic behavior of $K_{\mu}(z)$ as $z\rightarrow \infty$ (see e.g. \cite{71}) gives
\begin{proposition}When $\alpha=2$ and $\nu>0$, the asymptotic expansion of the function $q_{\alpha,\nu, n}(r)$ when $r\rightarrow \infty$ is given by
\begin{equation*}
q_{2,\nu, n}(r)\sim \frac{e^{-r}}{2^{\frac{n-1}{2}+\nu}\pi^{\frac{n-1}{2}}\Gamma(\nu)}\sum_{j=0}^{N-1} \frac{1}{2^j j!}\frac{\Gamma\left(\frac{n+1}{2}+j-\nu\right)}{\Gamma\left(\frac{n+1}{2}-j-\nu\right)}r^{\nu-\frac{n+1}{2}-j} + O\left(r^{\nu-\frac{n+1}{2}-N} e^{-r}\right).
\end{equation*} In particular, the large--$r$ leading term of $q_{2,\nu, n}(r)$ is given by
\begin{equation*}
q_{2,\nu, n}(r)\sim \frac{r^{\nu-\frac{n+1}{2}}e^{-r}}{2^{\frac{n-1}{2}+\nu}\pi^{\frac{n-1}{2}}\Gamma(\nu)}.
\end{equation*}
\end{proposition}

For $r\rightarrow 0$,  the explicit series representation of $K_{\mu}(z)$ about $z=0$ (see e.g. \cite{71})  gives:
\begin{proposition}\label{pro11}When $\alpha=2$ and $\nu>0$, \\
I. If $\frac{n}{2}-\nu$ is not an integer, then the   series expansion of $q_{\alpha, \nu, n}(r)$ about $r=0$ is given by
\begin{equation*}
q_{2, \nu, n}(r) = \frac{1}{2^{n}\pi^{\frac{n-2}{2}}\Gamma(\nu)\sin \left[\pi\left(\frac{n}{2}-\nu\right)\right]} \sum_{j=0}^{\infty}\left(\frac{(r/2)^{2j+2\nu-n}}{j!\Gamma\left(j+\nu-\frac{n-2}{2}\right)}
-\frac{(r/2)^{2j}}{j!\Gamma\left(j+\frac{n+2}{2}-\nu\right)} \right).
\end{equation*}
II. If $\frac{n}{2}-\nu=l$ is an integer, then
\begin{equation*}\begin{split}
q_{2, \nu, n}(r) =& \frac{1}{2^{n}\pi^{\frac{n}{2}}\Gamma\left(\frac{n}{2}-l\right)}\Biggl( \sum_{j=0}^{|l|-1}(-1)^j \frac{(|l|-j-1)!}{j! }\left(\frac{r}{2}\right)^{2j-|l|-l}
\\&+(-1)^{l+1}\sum_{j=0}^{\infty}\frac{(r/2)^{|l|-l+2j}}{j!   (|l|+j)!}\left[ \log\left(\frac{r}{2}\right)^2- \psi(j+1)- \psi(|l|+j+1)\right] \Biggr).
\end{split}
\end{equation*}

\end{proposition}
Notice that the behavior of $q_{2, \nu, n}(r)$ as $r\rightarrow 0$ depends on whether $\frac{n}{2}-\nu$ is an integer. If $\frac{n}{2}-\nu$ is not an integer, $q_{2, \nu, n}(r)$ can be represented by the sum of two series, one is an absolutely convergent power series, and the other is the multiplication  of $r^{2\nu-n}$ with an absolutely convergent power series. When $\frac{n}{2}-\nu$ is an integer,  $q_{2, \nu, n}(r)$ can be represented by the sum of three terms, one is an absolutely convergent power series, one is the multiplication of $r^{\min\{2\nu-n, 0\}}$ with a polynomial, and the third one is the multiplication of $\ln r$ with an absolutely convergent power series.
\begin{corollary}\label{co1}When $\alpha=2$ and $\nu>0$, the $r\rightarrow 0$ leading term of the function $q_{\alpha, \nu, n}(r)$ depends on the value of $\nu$:\\
\begin{equation*}\begin{split} & \text{I. If}\;   \nu <\frac{n}{2},\;\; \text{then}\;\;
q_{2, \nu, n}\sim \frac{\Gamma\left(\frac{n}{2}-\nu\right)}{2^{2\nu}\pi^{\frac{n}{2}}\Gamma(\nu)}r^{2\nu -n}.\\
& \text{ II. If}\;   \nu > \frac{n}{2}, \;\;\text{then}\;\;
q_{2, \nu, n}\sim \frac{\Gamma\left(\nu-\frac{n}{2}\right)}{2^n\pi^{\frac{n}{2}}\Gamma(\nu)}.\\
& \text{III. If}\;  \nu =\frac{n}{2}, \;\;\text{then}\;\; q_{2, \nu, n}\sim \frac{1}{2^{n}\pi^{\frac{n}{2}}\Gamma\left(\frac{n}{2}\right)}\left\{-\log  \left(\frac{r}{2}\right)^2+2\psi(1)\right\}.\hspace{4cm}\end{split}\end{equation*}
\end{corollary}

\subsection{$\alpha\in (0,2)$}

When $\alpha\in (0,2)$, the asymptotic expansion of $q_{\alpha, \nu, n}(r)$ as $r\rightarrow \infty$ can be obtained   from the representation of $q_{\alpha, \nu, n}(r)$ given by \eqref{eq10_29_2}.
\begin{proposition}\label{pro5}
When $\alpha\in (0, 2)$ and $\nu>0$, the following formula is valid: For $r\rightarrow\infty$,
\begin{equation}\label{eq10_30_3}\begin{split}
q_{\alpha,\nu, n}(r) \sim & \frac{1}{ \pi^{\frac{n+2}{2}}}\sum_{j=1}^{N-1}\frac{(-1)^{j-1}}{j!}\frac{\Gamma(\nu+j)}{\Gamma(\nu)}\sin\frac{\pi\alpha j}{2}2^{ \alpha j}\Gamma\left(\frac{ n+\alpha j}{2}\right)\Gamma\left(\frac{ 2 +\alpha j}{2}\right)r^{-\alpha j-n}\\
&+\mathcal{R}_N(r),\end{split}
\end{equation}where
\begin{equation}
\left| \mathcal{R}_N(r)\right|\leq \frac{2^{\alpha N}\Gamma(N+\nu)\Gamma\left(\frac{n+\alpha N}{2}\right)\Gamma\left(\frac{2+\alpha N}{2}\right)}{\pi^{\frac{n+2}{2}}N!\Gamma(\nu)\sin^{\nu+N}\frac{\pi\alpha}{2}}r^{-\alpha N-n}.
\end{equation}
In particular, the large--$r$ leading term of $q_{\alpha,\nu, n}(r)$ is given by
\begin{equation*}
q_{\alpha,\nu, n}(r)\sim  \frac{2^{\alpha} \nu}{ \pi^{\frac{n+2}{2}}}\sin\frac{\pi\alpha }{2}\Gamma\left(\frac{ n+\alpha }{2}\right)\Gamma\left(\frac{ 2 +\alpha }{2}\right)r^{-\alpha -n}.
\end{equation*}
\end{proposition}
\begin{proof}

Making a change of variable in \eqref{eq10_29_2}, we have
\begin{equation}\label{eq10_30_1}\begin{split}
q_{\alpha,\nu, n}(r) =&-\frac{r^{-n}}{2^{\frac{n-2}{2}}\pi^{\frac{n+2}{2}}}\text{Im}\,\int_0^{\infty}
\frac{K_{\frac{n-2}{2}}(
u)}{\left(1+e^{\frac{i\pi\alpha}{2}}r^{-\alpha}u^{\alpha}\right)^{\nu}}u^{\frac{n}{2}}du.\end{split}
\end{equation}
Using integration by parts, one can check easily that for an infinitely differentiable function $g(y)$, we have
\begin{equation*}
g(y)=\sum_{j=0}^{N-1} \frac{g^{j}(0)}{j!} y^j + \frac{1}{(N-1)!}\int_0^y (y-v)^{N-1}g^{(N)}(v)dv.
\end{equation*}This shows that
\begin{equation}\label{eq10_30_2}
\frac{1}{\left(1+e^{\frac{i\pi\alpha}{2}}r^{-\alpha}u^{\alpha}\right)^{\nu}}=\sum_{j=0}^{N-1}\frac{(-1)^j}{j!}\frac{\Gamma(\nu+j)}{\Gamma(\nu)}e^{\frac{i\pi\alpha j}{2}}r^{-\alpha j}u^{\alpha j}+R_N(r, u),
\end{equation}where
\begin{equation*}\begin{split}
R_N(r,u) = & \frac{(-1)^N\Gamma(N+\nu)}{(N-1)!\Gamma(\nu)}\int_0^{e^{\frac{i\pi\alpha}{2}}r^{-\alpha}u^{\alpha}}\frac{\left(e^{\frac{i\pi\alpha}{2}}r^{-\alpha}u^{\alpha}-v\right)^{N-1}}
{\left(1+v\right)^{\nu+N}}dv\\
=&\frac{(-1)^N\Gamma(N+\nu)e^{\frac{i\pi N\alpha}{2}}}{(N-1)!\Gamma(\nu)}\int_0^{r^{-\alpha}u^{\alpha}}\frac{\left( r^{-\alpha}u^{\alpha}-v\right)^{N-1}}
{\left(1+e^{\frac{i\pi\alpha}{2}}v\right)^{\nu+N}}dv.\end{split}
\end{equation*}Notice that
\begin{equation*}
\left\vert\frac{1}{\left(1+e^{\frac{i\pi\alpha}{2}}v\right)^{\nu+N}}\right\vert \leq \frac{1}{\sin^{\nu+N}\frac{\pi\alpha}{2}}.
\end{equation*}Therefore
\begin{equation}\label{eq10_30_4}
\left\vert R_N(r,u) \right\vert\leq \frac{ \Gamma(N+\nu) }{N!\Gamma(\nu)\sin^{\nu+N}\frac{\pi\alpha}{2}}r^{-\alpha N}u^{\alpha N}.
\end{equation}Substituting \eqref{eq10_30_2} into \eqref{eq10_30_1}, we find that
\begin{equation*}\begin{split}
q_{\alpha,\nu, n}(r) =&-\frac{r^{-n}}{2^{\frac{n-2}{2}}\pi^{\frac{n+2}{2}}}\text{Im}\Biggl\{\sum_{j=0}^{N-1} \frac{(-1)^j}{j!}\frac{\Gamma(\nu+j)}{\Gamma(\nu)}e^{\frac{i\pi\alpha j}{2}}r^{-\alpha j} \int_0^{\infty} K_{\frac{n-2}{2}}(
u)u^{\frac{n}{2}+\alpha j}du\\&+\int_0^{\infty}K_{\frac{n-2}{2}}(
u)R_{N}(r, u) u^{\frac{n}{2}} du  \Biggr\}.
\end{split}\end{equation*}The assertion \eqref{eq10_30_3} follows from the formula (see e.g. \cite{71})
\begin{equation}\label{eq10_31_1}
\int_0^{\infty} K_{\frac{n-2}{2}}(
u)u^{\frac{n}{2}+\alpha j}du = 2^{\frac{n}{2}+\alpha j-1}\Gamma\left(\frac{ n+\alpha j}{2}\right)\Gamma\left(\frac{ 2 +\alpha j}{2}\right)
\end{equation}and the bound \eqref{eq10_30_4}.
\end{proof}
An immediate corollary of Proposition
\ref{pro5} is
\begin{corollary}
When $\alpha\in (0, 2)$ and $\nu>0$, the asymptotic expansion of $q_{\alpha,\nu, n}(r)$ as $r\rightarrow \infty$ is given by
\begin{equation}\label{eq10_30_3_1}\begin{split}
q_{\alpha,\nu, n}(r) \sim & \frac{1}{ \pi^{\frac{n+2}{2}}}\sum_{j=1}^{\infty}\frac{(-1)^{j-1}}{j!}\frac{\Gamma(\nu+j)}{\Gamma(\nu)}\sin\frac{\pi\alpha j}{2}2^{ \alpha j}\Gamma\left(\frac{ n+\alpha j}{2}\right)\Gamma\left(\frac{ 2 +\alpha j}{2}\right)r^{-\alpha j-n}.\end{split}
\end{equation}
\end{corollary}

For the special cases where $\nu=1$ or $n=1$, \eqref{eq10_30_3_1} was obtained in \cite{10, 18, 20, 19, 17} by different methods. Notice that the large $r$-asymptotic behavior of $q_{\alpha, \nu, n}(r)$ is very different for the case of $\alpha\in (0,2)$ and $\alpha=2$.    $q_{\alpha, \nu, n}(r)$  decays polynomially when $\alpha\in (0,2)$ and exponentially when $\alpha=2$. Also notice that the large-$r$ leading term of $q_{\alpha,\nu, n}(r)$ is of order $r^{-\alpha-n}$, which is the same as the large-$r$ leading term of $\tilde{s}_{\alpha,n}(r)$ \eqref{eq11_12_2}.

Now we turn to the asymptotic behavior of $q_{\alpha, \nu, n}(r)$ when $r\rightarrow 0$. For this purpose, we need to derive an integral representation of $q_{\alpha,\nu, n}(r)$ which extends the results of \cite{10, 18, 20, 19, 17}. We begin with the following lemmas.
\begin{lemma}\label{lem1}
Let $z$ be a complex number such that $\left|\text{arg}\, (z)\right|< \pi$, then
\begin{equation}\label{eq10_30_7}
\frac{1}{(1+ z)^{\nu}}=\frac{1}{\Gamma(\nu)}\int_{c-i\infty}^{c+i\infty} \Gamma(w)\Gamma(\nu-w)z^{-w}dw,
\end{equation}where $0<c<\nu$.
\end{lemma}\begin{proof}Using the fact that when $u\gg 1$ (see e.g. \cite{73}),
\begin{equation}\label{eq11_04_1}
|\Gamma(c+iu)|=\sqrt{2\pi} |u|^{c-\frac{1}{2}}e^{-\frac{\pi |u|}{2}}\left(1+O\left(\frac{1}{|u|}\right)\right),
\end{equation}we find that if $u\gg 1$,
\begin{equation}\label{eq11_12_3}
\Bigl|\Gamma(c+iu)\Gamma(\nu-c-iu)z^{-c-iu}\Bigr|\leq C|z|^{-c} e^{|u||\text{arg}\,z|}|u|^{\nu-1}e^{-\pi |u|}.
\end{equation}Here and the followings, $C$ or $C_1, C_2$ represent constants whose values can be different in different lines. \eqref{eq11_12_3} shows that the right hand side of \eqref{eq10_30_7} indeed defines an analytic function of $z$ when $\left|\text{arg}\, (z)\right|< \pi$. Now using
\begin{equation*}
e^{-z} =\frac{1}{2\pi i}\int_{c-i\infty}^{c+i\infty} \Gamma(w) z^{-w} dw,
\end{equation*}we have
\begin{equation}\label{eq10_30_8}\begin{split}
\frac{1}{(1+ z)^{\nu}}=&\frac{1}{\Gamma(\nu)}\int_0^{\infty} y^{\nu-1}e^{-y(1+z)}dy\\
=&\frac{1}{\Gamma(\nu)}\int_0^{\infty} y^{\nu-1}e^{-y}\left\{\frac{1}{2\pi i}\int_{c-i\infty}^{c+i\infty}  \Gamma(w) y^{-w} z^{-w}dw\right\}dy.
\end{split}
\end{equation}For $u\gg 1$, \eqref{eq11_04_1} implies that
\begin{equation*}
\Bigl\vert e^{-y} y^{\nu-c-iu-1}\Gamma(c+iu) z^{-c-iu}\Bigr\vert \leq C e^{|u||\text{arg}\,z|} |z|^{-c} e^{-y} y^{\nu -c-1} e^{-\frac{\pi}{2}|u|}|u|^{c-\frac{1}{2}}.
\end{equation*}Therefore we can interchange the order of integrations in \eqref{eq10_30_8} when $\left|\text{arg}\, (z)\right|< \pi/2$ to obtain \eqref{eq10_30_7}. The general case where $\left|\text{arg}\, (z)\right|< \pi$ follows by analytic continuation.
\end{proof}
\begin{lemma}\label{lem2}For $\alpha\in (0,2)$ and $u>0$, we have
\begin{equation*}
\text{Im}\; \frac{-1}{\left(1+ e^{\frac{i\pi\alpha}{2}}u^{\alpha}\right)^{\nu}}=\frac{1}{2\pi i} \frac{1}{\Gamma(\nu)}\int_{c-i\infty}^{c+i\infty} \Gamma(w)\Gamma(\nu-w) u^{-\alpha w} \sin\frac{\pi\alpha w}{2} dw,
\end{equation*}where $0<c<\nu$.

\end{lemma}\begin{proof}
From Lemma \ref{lem1}, we have
\begin{equation*}\begin{split}
\text{Im}\; \frac{-1}{\left(1+ e^{\frac{i\pi\alpha}{2}}u^{\alpha}\right)^{\nu}}=&\text{Im}\,\left\{\frac{-1}{2\pi i}\frac{1}{\Gamma(\nu)}\int_{c-i\infty}^{c+i\infty}\Gamma(w)\Gamma(\nu-w) e^{-\frac{i\pi \alpha w}{2}}u^{-\alpha w}dw\right\}\\
=&\text{Im}\,\left\{\frac{-1}{2\pi i}\frac{1}{\Gamma(\nu)}\int_{c-i\infty}^{c+i\infty}\Gamma(w)\Gamma(\nu-w) \cos\left(\frac{\pi\alpha w}{2}\right)u^{-\alpha w}dw\right\}\\
&+\text{Re}\,\left\{\frac{1}{2\pi i}\frac{1}{\Gamma(\nu)}\int_{c-i\infty}^{c+i\infty}\Gamma(w)\Gamma(\nu-w) \sin\left(\frac{\pi \alpha w}{2}\right)u^{-\alpha w}dw\right\}.\end{split}
\end{equation*}
Therefore, it suffices to show that
\begin{equation}\label{eq10_31_2}
\text{Im}\,\left\{\frac{-1}{2\pi i}\frac{1}{\Gamma(\nu)}\int_{c-i\infty}^{c+i\infty}\Gamma(w)\Gamma(\nu-w) \cos\left(\frac{\pi\alpha w}{2}\right)u^{-\alpha w}dw\right\}=0
\end{equation}and
\begin{equation}\label{eq10_31_3}
\text{Im}\,\left\{\frac{-1}{2\pi i}\frac{1}{\Gamma(\nu)}\int_{c-i\infty}^{c+i\infty}\Gamma(w)\Gamma(\nu-w) \sin\left(\frac{\pi\alpha w}{2}\right)u^{-\alpha w}dw\right\}=0.
\end{equation}Notice that the functions $z\mapsto\sin z$, $z\mapsto\cos z$, $z\mapsto u^z, \,u\in \mathbb{R}$ and $z\mapsto \Gamma(z)$ are all real on the real axis. Reflection principle implies that they satisfy $f(\bar{z})=\overline{f(z)}$. On the other hand,
\begin{equation}\label{eq10_31_4}\begin{split}
&\frac{-1}{2\pi i}\frac{1}{\Gamma(\nu)}\int_{c-i\infty}^{c+i\infty}\Gamma(w)\Gamma(\nu-w) \cos\left(\frac{\pi\alpha w}{2}\right)u^{-\alpha w}dw\\=&-\frac{1}{2\pi}\int_{-\infty}^{\infty}\Gamma(c+iy)\Gamma(\nu-c-iy)\cos\left(\frac{\pi\alpha (c+iy)}{2}\right)u^{-\alpha (c+iy)}dy.
\end{split}\end{equation}By taking the complex conjugate of \eqref{eq10_31_4}, it is easy to see that \eqref{eq10_31_4} is real. Therefore  \eqref{eq10_31_2} holds. \eqref{eq10_31_3} follows analogously.

\end{proof}

Now we can prove the following.
\begin{proposition}
The function $q_{\alpha,\nu, n}(r)$ has the following integral representation:
\begin{equation}\label{eq10_30_6}
\begin{split}q_{\alpha, \nu, n}(r) = \frac{r^{-n}}{\pi^{\frac{n}{2}}\Gamma(\nu)}\frac{1}{2\pi i} \int_{c-iw}^{c+iw}f_{\alpha,\nu, n}(w; r)dw,
\end{split}
\end{equation}where $0 <c<\min\{\nu, \frac{1}{\alpha}\}$. Here
\begin{equation}\label{eq11_06_1}f_{\alpha,\nu, n}(w; r)= \Gamma(w)\Gamma(\nu-w) \left(\frac{r}{2}\right)^{\alpha w}\frac{ \Gamma\left(\frac{n-\alpha w}{2}\right)}{\Gamma\left(\frac{\alpha w}{2}\right)}.\end{equation}
\end{proposition}
\begin{proof}
Using \eqref{eq10_29_2} and Lemma \ref{lem2}, we have
\begin{equation}\label{eq10_30_9}
q_{\alpha,\nu, n}(r)=\frac{r^{\frac{2-n}{2}}}{2^{\frac{n-2}{2}}\pi^{\frac{n+2}{2}}} \int_0^{\infty}K_{\frac{n-2}{2}}(ru) u^{\frac{n}{2}}\left\{\frac{1}{2\pi i}\frac{1}{\Gamma(\nu)}\int_{c-i\infty}^{c+i\infty}\Gamma(w)\Gamma(\nu-w) \sin\frac{\pi \alpha w}{2}u^{-\alpha w} \right\}du.
\end{equation}For $0\leq u\leq 1$ and $y\gg 1$, \eqref{eq11_04_2} and \eqref{eq11_04_1} imply that
\begin{equation}\label{eq11_13_1}
\Bigl|K_{\frac{n-2}{2}}(ru) u^{\frac{n}{2}-\alpha c-i\alpha y}\Gamma(c+iy)\Gamma(\nu-c-iy) \sin\frac{\pi \alpha (c+iy)}{2}\Bigr|\leq C u^{-\alpha c} |y|^{\nu-1}e^{-\frac{\pi|y|}{2}\left(2-\alpha\right)}.
\end{equation}For $  u\geq 1$ and $y\gg 1$, \eqref{eq11_04_2} and \eqref{eq11_04_1} give
\begin{equation}\label{eq11_13_2}
\Bigl|K_{\frac{n-2}{2}}(ru) u^{\frac{n}{2}-\alpha c-i\alpha y}\Gamma(c+iy)\Gamma(\nu-c-iy) e^{-\frac{i\pi \alpha (c+iy)}{2}}\Bigr|\leq C u^{\frac{n}{2}-\alpha c} e^{-ru} |y|^{\nu-1}e^{-\frac{\pi|y|}{2}(2-\alpha)}.
\end{equation}Since $c<1/\alpha$, \eqref{eq11_13_1} and \eqref{eq11_13_2} show that we can interchange the order of integrations in \eqref{eq10_30_9} and use \eqref{eq10_31_1} to obtain
\begin{equation*}
q_{\alpha,\nu, n}(r)=\frac{r^{-n}}{\pi^{\frac{n+2}{2}}\Gamma(\nu)}\frac{1}{2\pi i}\int_{c-i\infty}^{c+i\infty}\Gamma(w)\Gamma(\nu-w) \sin\frac{\pi \alpha w}{2}\left(\frac{r}{2}\right)^{\alpha w}\Gamma\left(\frac{n-\alpha w}{2}\right)\Gamma\left(\frac{2-\alpha w}{2}\right)dw.
\end{equation*}
 \eqref{eq10_30_6} follows by using the identity
\begin{equation*}
\Gamma(z)\Gamma(1-z)=\frac{\pi}{\sin \pi z}.
\end{equation*}
\end{proof}By making a change of variable, it is easy to see that \eqref{eq10_30_6} coincides with the results of \cite{10, 18, 20, 19, 17} in the special cases where $\nu=1$ or $n=1$. However, instead of proving that the right hand side of \eqref{eq10_30_6} satisfies the equation \eqref{eq10_29_4}, we choose to give a direct derivation of \eqref{eq10_30_6} here.

To find the asymptotic expansion of $q_{\alpha, \nu, n}(r)$ as $r\rightarrow 0$, we can apply residue calculus to \eqref{eq10_30_6}. More precisely, we shift the contour of integration to the right   and pick up the residues of  $f_{\alpha,\nu, n}(w; r)$.
Since the function $\Gamma(z)$ has simple poles at the point $z=-k$, $k=0, 1, 2, ...$ with residue
 \begin{equation}\label{eq11_13_3}\text{Res}_{z=-k}\Gamma(z)=\frac{(-1)^k}{k!},\end{equation}
  the poles of $f_{\alpha,\nu, n}(w; r)$   that are situated to the right of the line $\text{Re}\; w=c$  are at the points
$$w= \frac{n+2j}{\alpha}, \;\;j=0, 1, 2,... \;\;\;\text{and}\;\;\; w= \nu+l,\;\; l=0,1,2,....$$ These are simple poles except if $n+2j =\alpha(\nu+l)$ for some $(j, l)\in \hat{\mathbb{N}}\times\hat{\mathbb{N}}$. Here $\hat{\mathbb{N}}=\mathbb{N}\cup\{0\}$ is the set of nonnegative integers. Define the set $\Lambda_{n}$ by
\begin{equation*}
\Lambda_n =\left\{ (\alpha, \nu)\in (0,2)\times \mathbb{R}_+\,|\, n+2j =\alpha(\nu+l)\;\;\text{for some }\; (j, l) \in \hat{\mathbb{N}}\times\hat{\mathbb{N}}\right\},
\end{equation*}so that $(\alpha, \nu)\in \Lambda_n$ if and only if $f_{\alpha, \nu, n}(w; r)$ has double poles on the right half plane $\text{Re}\,w>0$. The detail  characterization of the set $\Lambda_n$ is deferred to the end of this section. In the following, we give the asymptotic expansion of $q_{\alpha,\nu, n}(r)$ when $r\rightarrow 0$.

\begin{proposition}\label{pro4}
When $\alpha\in (0,2)$ and $\nu>0$, the following formula is valid:
\begin{equation*}
q_{\alpha,\nu, n}(r)= S_{\alpha, \nu, n;N}(r)+R_{\alpha, \nu, n; N}(r),
\end{equation*}
where
\begin{equation}\label{eq10_31_10}\begin{split}
&S_{\alpha, \nu, n;N}(r)\\=&\frac{1}{2^n\pi^{\frac{n}{2}}\Gamma(\nu)}\Biggl\{ \sum_{\substack{0\leq l\leq N-1\\ l \neq \frac{n+2j}{\alpha}-\nu\;\\\text{for any}\;j\in\hat{\mathbb{N}}} }\frac{(-1)^{l}}{l!}\Gamma(\nu+l)\frac{\Gamma\left(\frac{n-\alpha(\nu+l)}{2}\right)}{\Gamma\left(\frac{\alpha(\nu+l)}{2}\right)} \left(\frac{r}{2}\right)^{\alpha \nu-n+\alpha l }\\&+\frac{2}{\alpha}\sum_{\substack{0\leq j\leq \left[\frac{\alpha\left(\nu+N-\frac{1}{2}\right)-n}{2}\right]\\j\neq \frac{\alpha(\nu+l)-n}{2}\\ \text{for any}\;l\in\hat{\mathbb{N}}}}\frac{(-1)^j}{j!}
\frac{\Gamma\left(\frac{n+2j}{\alpha}\right)}{\Gamma\left(\frac{n+2j}{2}\right)}\Gamma\left(\frac{\alpha \nu -n-2j}{\alpha}\right)\left(\frac{r}{2}\right)^{2j}\\
&- \sum_{\substack{l= \frac{n+2j}{\alpha}-\nu\;\\\text{for some}\;j\in\hat{\mathbb{N}}\\0\leq l\leq N-1}}\frac{(-1)^{l +\frac{\alpha(\nu+l)-n}{2}}}{l!\left(\frac{\alpha(\nu+l)-n}{2}\right)!}\frac{\Gamma(\nu+l)}{\Gamma\left(\frac{\alpha(\nu+l)}{2}\right)}
\Biggl\{\frac{2}{\alpha}\psi(\nu+l)-\frac{2}{\alpha}\psi(l+1)-
\psi\left(\frac{\alpha(\nu+l)}{2}\right)\\&- \psi\left(\frac{\alpha(\nu+l)-n}{2}+1\right)+ \log\left(\frac{r}{2}\right)^2\Biggr\}\left(\frac{r}{2}\right)^{\alpha\nu-n+\alpha l}\Biggr\},\end{split}
\end{equation} and
\begin{equation}\label{eq11_07_1}
R_{\alpha, \nu, n; N}(r)=O\left(r^{\alpha\nu-n+\alpha \left(N-\frac{3}{4}\right)}\right)\hspace{1cm}\text{as}\;\;r\rightarrow 0.
\end{equation}
Note that the third summation in \eqref{eq10_31_10} is zero if and only if $(\alpha,\nu)\notin \Lambda_n$.
\end{proposition}
\begin{proof}We use the integral representation of $q_{\alpha,\nu, n}(r)$ given by \eqref{eq10_30_6}.
Given a positive integer $N$,  we   observe that there is at most one point of the form $w=\frac{n+2j}{\alpha}$, $j\in \hat{\mathbb{N}}$ in between $\nu+N-1$ and $\nu+N$ since $\alpha<2$. If there is no such point, let $Q_N =\nu + N-\frac{1}{2}$. If $w_N =\frac{n+2j_N}{\alpha}$, $j_N\in \hat{\mathbb{N}}$, lies between $\nu+N-1$ and $\nu + N-\frac{1}{2}$, let $Q_N =\nu+N-\frac{1}{4}$. If $w_N =\frac{n+2j_N}{\alpha}$, $j_N\in \hat{\mathbb{N}}$, lies between $\nu + N-\frac{1}{2}$ and $\nu +N$, let $Q_N =\nu+N-\frac{3}{4}$. Notice that \begin{equation}\label{eq11_07_12_1}\min\left\{\left|Q_N -(\nu+N-1)\right|, \left|Q_N-(\nu+N)\right| \right\}\geq \frac{1}{4},\end{equation} and if there exists $w_N$ lying between $\nu+N-1$ and $\nu+N$, then
 \begin{equation}
 \label{eq11_07_12_2}\min\left\{  \left|Q_N-w_N\right|\right\}\geq \frac{1}{4}.
\end{equation} Consider the contour $\mathcal{C}_{N, M}$ which consists of the boundary of the rectangle \\$\left\{w\,|\, c\leq  \text{Re}\, w \leq Q_N, \;\; |\text{Im}\,w|\leq M\right\}$. The poles of $f_{\alpha,\nu, n}(w; r)$ inside $\mathcal{C}_{N, M}$ are the points $w=\nu+l$, $l=0, 1, \ldots, N-1$ and $w=\frac{n+2j}{\alpha}$, $j=0, 1,\ldots, \left[\frac{\alpha(\nu+N-\frac{1}{2})-n}{2}\right]$.  When $c\leq  \text{Re}\, w \leq  Q_N$ and $\text{Im}\, w=\pm M$,  \eqref{eq11_04_1} shows that for $M\gg 1$,
\begin{equation*}
\left| f_{\alpha,\nu,n}(w;r)\right|\leq C M^{\frac{n}{2}-\alpha c+\nu-1}e^{-\pi M} 2^{-\alpha c} \max\{1, r\}^{\alpha(\nu+N)}.
\end{equation*}Therefore, the integral of $f_{\alpha, \nu, n}(w; r)$ over the two horizontal segments of the rectangle $\mathcal{C}_{N, M}$ goes to zero as $M\rightarrow \infty$. This gives
\begin{equation*}\begin{split}
\frac{1}{2\pi i}\int_{c-i\infty}^{c+i\infty} f_{\alpha,\nu, n}(w; r) dw=& -\sum_{z\;\text{inside}\,\mathcal{C}_{N, M}} \text{Res}_{w=z}f_{\alpha,\nu, n}(w; r)+\frac{1}{2\pi i}\int_{ Q_N-i\infty}^{ Q_N+i\infty}f_{\alpha,\nu, n}(w,r)dw.\end{split}
\end{equation*}
If $(\alpha, \nu)\notin \Lambda_n$, all the poles of $f_{\alpha,\nu, n}(w; r) $ inside $\mathcal{C}_{N, M}$ are simple and \eqref{eq11_13_3} implies that
\begin{equation}\label{eq10_31_11}\begin{split}
&-\sum_{z\;\text{inside}\,\mathcal{C}_{N, M}} \text{Res}_{w=z}f_{\alpha,\nu, n}(w; r)\\=&-\sum_{l=0}^{N-1}\text{Res}_{w=\nu+l}f_{\alpha,\nu, n}(w; r)-\sum_{j=0}^{\left[\frac{\alpha\left(\nu+N-\frac{1}{2}\right)-n}{2}\right]}\text{Res}_{w=\frac{n+2j}{\alpha}}f_{\alpha,\nu, n}(w; r)\\
=&\sum_{l=0}^{N-1}\frac{(-1)^{l}}{l!}\Gamma(\nu+l)\frac{\Gamma\left(\frac{n-\alpha(\nu+l)}{2}\right)}{\Gamma\left(\frac{\alpha(\nu+l)}{2}\right)} \left(\frac{r}{2}\right)^{\alpha(\nu+l)}\\&+\frac{2}{\alpha}\sum_{j=0}^{\left[\frac{\alpha\left(\nu+N-\frac{1}{2}\right)-n}{2}\right]}\frac{(-1)^j}{j!}
\frac{\Gamma\left(\frac{n+2j}{\alpha}\right)}{\Gamma\left(\frac{n+2j}{2}\right)}\Gamma\left(\frac{\alpha \nu -n-2j}{\alpha}\right)\left(\frac{r}{2}\right)^{n+2j}.
\end{split}\end{equation}This gives \eqref{eq10_31_10} when $(\alpha,\nu)\notin\Lambda_n$. On the other hand, if $(\alpha, \nu)\in \Lambda_n$, then there are $(j, l)\in \hat{\mathbb{N}}\times\hat{\mathbb{N}}$ such that $n+2j =\alpha(\nu+l)$. In other words, $l=\frac{n+2j}{\alpha}-\nu$ for some $j\in \hat{\mathbb{N}}$ and $j=\frac{\alpha(\nu+l)-n}{2}$ for some $l\in \hat{\mathbb{N}}$. We have to omit such terms from the first and second summations in \eqref{eq10_31_11} and add the residues of $f_{\alpha, \nu, n}(w;r)$ at $w= \nu+l=\frac{n+2j}{\alpha}$, which gives
\begin{equation}\begin{split}\label{eq11_13_4}
-\frac{2}{\alpha}\sum_{\substack{l= \frac{n+2j}{\alpha}-\nu\;\\\text{for some}\;j\in\hat{\mathbb{N}}\\0\leq l\leq N-1}}\frac{(-1)^{l +\frac{\alpha(\nu+l)-n}{2}}}{l!\left(\frac{\alpha(\nu+l)-n}{2}\right)!}&\frac{\Gamma(\nu+l)}{\Gamma\left(\frac{\alpha(\nu+l)}{2}\right)}\Biggl\{\psi(\nu+l)-\psi(l+1)-\frac{\alpha}{2}
\psi\left(\frac{\alpha(\nu+l)}{2}\right)\\&-\frac{\alpha}{2}\psi\left(\frac{\alpha(\nu+l)-n}{2}+1\right)+\alpha\log\frac{r}{2}\Biggr\}\left(\frac{r}{2}\right)^{\alpha(\nu+l)}.
\end{split}\end{equation}This gives the third term in \eqref{eq10_31_10} when $(\alpha,\nu)\in\Lambda_n$.  To obtain  \eqref{eq11_13_4}, we  have used the formula
\begin{equation*}
\Gamma(-k + z) = \frac{(-1)^k}{k!}\left\{ \frac{1}{z}+\psi(k+1)\right\}+O(z)\;\;\text{as}\;z\rightarrow 0,\hspace{1cm}k\in\hat{\mathbb{N}}.
\end{equation*}  Finally we consider the remainder term.
For   $y$ large enough, \eqref{eq11_04_1} gives
\begin{equation}\label{eq11_06_2}
\Bigl|f_{\alpha, \nu, n}(Q_N+iy; r) \Bigr|\leq C 2^{-\alpha Q_N}r^{\alpha Q_N} e^{-\pi |y|}|y|^{\frac{n}{2}-\alpha Q_N+\nu -1}.
\end{equation}Therefore
\begin{equation*}
\left|\frac{1}{2\pi i}\int_{ Q_N-i\infty}^{ Q_N+i\infty}f_{\alpha,\nu, n}(w,r)dw\right|\leq Cr^{\alpha Q_N},
\end{equation*}which gives \eqref{eq11_07_1}.

\end{proof}By setting $n=1$ or $\nu=1$, one can check that the results of \cite{20, 19} and \cite{18} were recovered respectively. From \eqref{eq10_31_10}, we see that $q_{\alpha,\nu, n}(r)\rightarrow \infty$ as $r\rightarrow 0$ if and only if $\alpha\nu\leq n$. More precisely, we have
\begin{corollary}\label{co2}
Let $\alpha\in (0,2)$ and $\nu>0$. As $r\rightarrow 0$, the leading term of $q_{\alpha,\nu, n}(r)$ is given by
\begin{equation*}
\begin{split}
&\text{I.} \;\;\text{If}\;\alpha\nu <n, \;\;\text{then}\;\; q_{\alpha,\nu,n}(r)\sim \frac{1}{2^{\alpha\nu}\pi^{\frac{n}{2}} }\frac{\Gamma\left(\frac{n-\alpha\nu}{2}\right)}{\Gamma\left(\frac{\alpha\nu}{2}\right)} r^{\alpha\nu-n}\\
&\text{II.} \;\;\text{If}\;\alpha\nu >n, \;\;\text{then}\;\; q_{\alpha,\nu,n}(r)\sim \frac{1}{2^{n-1}\alpha\pi^{\frac{n}{2}}\Gamma(\nu)}\frac{\Gamma\left(\frac{n}{\alpha}\right)}{\Gamma\left(\frac{n}{2}\right)}\Gamma\left(\frac{\alpha\nu-n}{\alpha}\right).\\
&\text{III.} \;\;\text{If}\;\alpha\nu =n, \;\;\text{then}\;\; q_{\alpha,\nu,n}(r)\sim \frac{1}{2^{n}\pi^{\frac{n}{2}}\Gamma\left(\frac{n}{2}\right)}\left\{-\log \left(\frac{r}{2}\right)^2+\psi\left(\frac{n}{2}\right)+\psi(1)-\frac{2}{\alpha}\left(\psi(\nu)-\psi(1)\right)\right\}.
\end{split}
\end{equation*}
\end{corollary}One notice that putting $\alpha = 2$ in the results of Corollary \ref{co2}, one obtains the results of Corollary \ref{co1}. We also remark that the large-$r$ asymptotic expansion of $q_{\alpha, \nu, n}(r)$ when $\alpha\in (0,2)$  \eqref{eq10_30_3} can also be obtained from \eqref{eq10_30_6} by considering the poles on the left of the line $\text{Re}\; w=c$.

Now we analyze in more detail the conditions for which $(\alpha, \nu)\in \Lambda_n$. Observe that it is necessary that both $\alpha$ and $\nu$ are  rational numbers or both are irrational numbers. If $\alpha$ and $\nu$ are both irrational, then
\begin{equation*}
n+2j_1=\alpha(\nu+l_1) \hspace{0.5cm}\text{and}\hspace{0.5cm}n+2j_2=\alpha(\nu+l_2)\hspace{1cm}\text{for}\;\; (j_1, l_1), (j_2, l_2)\in \hat{\mathbb{N}}\times \hat{\mathbb{N}}
\end{equation*}if and only if $$2(j_1-j_2)=\alpha(l_1-l_2),$$if and only if $j_1=j_2$ and $l_1=l_2$.  In other words, if $\alpha$ and $\nu$ are both irrational numbers, there is at most one pair of nonnegative integers $(j, l)$ such that $n + 2j =\alpha(\nu+l)$. For the case where both $\alpha$ and $\nu$ are rational numbers, we write $\alpha = a/b$ and $\nu = e/f$, where $a, b, e, f\in \mathbb{N}$, $\text{gcd}\,(a, b)=1$ and $\text{gcd}\,(e,f)=1$. Then
\begin{equation}\label{eq11_04_6}
n+2j = \alpha(\nu +l) \Longrightarrow n+2j =\frac{a}{b}\frac{(e+lf)}{f}, \hspace{1cm}(j, l)\in \mathbb{N}\times\hat{\mathbb{N}}.
\end{equation}Since $\text{gcd}\;(a,b)=1$, $\text{gcd}\,(e,f)=1$ and the left hand side is an integer, it is necessary that $f$ divides $a$ and $b$ divides $e+lf$. Let $a=mf$ for some $m\in \mathbb{N}$. Since $\text{gcd}\, (a,b)=1$, this implies that $\text{gcd}\,(b,f)=1$. Therefore, it is always possible to find   positive integers $k$ and $l$ such that \begin{equation}\label{eq11_04_5}e+lf = kb.\end{equation} We have then $$\frac{a}{b}\frac{(e+lf)}{f}=km.$$ However, there does not necessary exist a nonnegative integer $j$ such that $n+2j=km$. We need to discuss the parity of $n$. If $n$ is odd, then $n+2j$ is odd for any $j\in \hat{\mathbb{N}}$. Therefore for $n+2j=km$, it is necessary that $m$ and $k$ are both odd. $m$ is odd implies that both $a$ and $f$ are even or both are odd. In the case where $a$ and $f$ are both even, $b$ and $e$ are both odd. Then for any $l\in \hat{\mathbb{N}}$, $e+lf$ is always odd. Therefore  if $(k, l)\in \mathbb{N}\times\hat{\mathbb{N}}$ is a solution to \eqref{eq11_04_5}, $k$ is necessary odd. This implies that there are always solutions $(j, l)\in \hat{\mathbb{N}}\times \hat{\mathbb{N}}$ to \eqref{eq11_04_6}. Moreover, if $(j_0, l_0)$ is the smallest nonnegative solution, then all the other solutions are given by $(j_q, l_q)= (j_0+aq/2, l_0+bq), q\in \hat{\mathbb{N}}$. In the case where $a$ and $f$ are both odd,  if $(k, l)$ is a solution to \eqref{eq11_04_5}, then so is $(k+f, l+b)$, but one of $k$ and $k+f$ must be odd. In this situation, we also find that there always exist solutions to \eqref{eq11_04_6}. Moreover, if $(j_0, l_0)$ is the smallest nonnegative solution, then all the other solutions are given by $(j_q, l_q) =(j_0+aq, l_0+2bq), q\in \hat{\mathbb{N}}$.

For the case where $n$ is even, $n+2j$ is even. In this case, $m$ can be either odd or even. If $m$ is odd, $k$ has to be even. The same reasoning above shows that $a$ and $f$ has to be both odd and   there is always solutions to \eqref{eq11_04_6}. Moreover, if $(j_0, l_0)$ is the smallest nonnegative solution, then all other solutions are given by $(j_q, l_q)=(j_0+aq, l_0+2bq), q\in \hat{\mathbb{N}}$. If $m$ is even, then $a$ is even and there is no restriction on $k$. In this case, there is always solutions to \eqref{eq11_04_6} and all the solutions can be expressed by the smallest nonnegative solution $(j_0, l_0)$ via $(j_q, l_q)=(j_0+aq/2, l_0+bq), q\in \hat{\mathbb{N}}$.

We summarize the results as follows:
\begin{proposition}\label{pro15}
Let $\alpha \in (0,2)$ and $\nu >0$.  \\
I. If $\alpha$ and $\nu$ are both irrational numbers such that $(\alpha,\nu)\in \Lambda_n$, then there is a unique solution $(j, l)$ to the equation $n+2j=\alpha(\nu+l), \;\;(j,l)\in \hat{\mathbb{N}}\times \hat{\mathbb{N}}$.\\
II. If $\alpha$ and $\nu$ are rational numbers with $\alpha=a/b$, $\nu=e/f$, $a,b, e,f\in \mathbb{N}$ and $\text{gcd}\,(a,b)=\text{gcd}\,(e,f)=1$, then $(\alpha,\nu)\in \Lambda_n$ if and only if $f$ divides $a$ and one of the following holds:\\
\;\; (i) $n$ is odd and $a/f$ is odd.\\
\;\; (ii) $n$ is even and if $a/f$ is odd, then $a$ and $f$ are both odd. \\
Moreover, if $(j_0, l_0)$ is the smallest nonnegative solution of $n+2j=\alpha(\nu+l), \;\;(j,l)\in \hat{\mathbb{N}}\times \hat{\mathbb{N}}$, then if $a$ is odd, all   solutions are given by $(j_q, l_q)=(j_0+aq, l_0+2bq), q\in \hat{\mathbb{N}}$. If $a$ is even, all solutions are given by $(j_q, l_q)=(j_0+aq/2, l_0+bq), q\in \hat{\mathbb{N}}$.

\end{proposition}

\section{Representation of $q_{\alpha,\nu, n}(r)$ in terms of entire functions} In this section, we are going to analyze under what conditions the asymptotic expansion of $q_{\alpha, \nu, n}(r)$ given in Proposition \ref{pro4} will give rise to representation of $q_{\alpha,\nu, n}(r)$ in terms of absolutely convergent power series.
First we have

\begin{proposition}\label{pro6}
Let $\alpha  \in (0,2)$ and $\nu>0$. Then
 \begin{equation*}
 q_{\alpha, \nu, n}(r)=\lim_{N\rightarrow \infty} S_{\alpha,\nu, n; N}(r)
 \end{equation*}
uniformly for $r$ in any compact subsets of $(0, \infty)$. Here $S_{\alpha,\nu, n; N}(r)$ is given by \eqref{eq10_31_10}.

\end{proposition}
\begin{proof}
Following from the proof of Proposition \ref{pro4}, the remainder term $R_{\alpha,\nu,n;N}(r)= q_{\alpha,\nu,n}(r)- S_{\alpha,\nu, n; N}(r)$ is given by
\begin{equation*}\begin{split}
R_{\alpha,\nu,n;N}(r)=\frac{r^{-n}}{\pi^{\frac{n}{2}}\Gamma(\nu)} \frac{1}{2\pi i}\int_{Q_N-i\infty}^{Q_N+i\infty}f_{\alpha,\nu,n}(w;r)dw=\frac{r^{-n}}{\pi^{\frac{n}{2}}\Gamma(\nu)} \frac{1}{2\pi }\int_{-\infty}^{ \infty}f_{\alpha,\nu,n}(Q_N+iy ;r)dy,
\end{split}\end{equation*}where $f_{\alpha,\nu, n}(w;r)$ is given by \eqref{eq11_06_1}. We rewrite it as
\begin{equation}\label{eq11_06_3}\begin{split}
 f_{\alpha, \nu, n}(Q_N+iy; r) =& \frac{(r/2)^{\alpha(Q_N+iy)}\Gamma(Q_N +iy)}{\Gamma(1-\nu+Q_N+iy)
\Gamma\left(\frac{2+\alpha(Q_N+iy)-n}{2}\right)\Gamma\left(\frac{\alpha \left(Q_N+iy\right)}{2}\right)}\\&\times \frac{\pi^2}{\sin \pi \left(Q_N+iy -\nu\right)\sin   \frac{\pi\left( \alpha(Q_N+iy)-n\right)}{2}}.
\end{split}\end{equation}We want to find an upper bound for $| f_{\alpha, \nu, n}(Q_N+iy; r) |$ when $Q_N$ is large. Since $\overline{f_{\alpha, \nu, n}(Q_N+iy; r) }=f_{\alpha, \nu, n}(Q_N-iy; r)$, it suffices to consider $y\geq 0$. By Stirling's formula (see e.g. \cite{73}),
\begin{equation}\label{eq11_13_5}
\Gamma(z) \sim \sqrt{2\pi} z^{z-\frac{1}{2}}e^{-z}\hspace{0.3cm}\text{as}\;\; |z|\rightarrow \infty, \hspace{1cm}| \text{arg}\,z|<\pi.
\end{equation}Therefore, if $0\leq v\leq   u$,
\begin{equation*}
C_1 u^{u-\frac{1}{2}}e^{-u-\frac{\pi}{4} v} \leq |\Gamma(u+iv)|\leq C_2(\sqrt{2}u)^{u-\frac{1}{2}}e^{-u};
\end{equation*}whereas if $v\geq u>0$,
\begin{equation*}
C_1 u^{u-\frac{1}{2}}e^{-u-\frac{\pi}{2} v} \leq |\Gamma(u+iv)|\leq C_2(\sqrt{2}v )^{u-\frac{1}{2}}e^{-u-\frac{\pi}{4}v}.
\end{equation*}These imply that if $Q_N\gg 1$ and  $0\leq y \leq Q_N$,
\begin{equation*}
\left| \frac{ \Gamma(Q_N +iy)}{\Gamma(1-\nu+Q_N+iy)
\Gamma\left(\frac{2+\alpha(Q_N+iy)-n}{2}\right)\Gamma\left(\frac{\alpha \left(Q_N+iy\right)}{2}\right)}\right|\leq C\frac{\left(\sqrt{2}(2e/\alpha)^{\alpha}e^{\left(\alpha+1\right)\frac{\pi}{4}}\right)^{Q_N}}{Q_N^{\alpha Q_N+1-\nu-\frac{n}{2}}};
\end{equation*}and if $ Q_N\leq y$,
\begin{equation*}
\left| \frac{ \Gamma(Q_N +iy)}{\Gamma(1-\nu+Q_N+iy)
\Gamma\left(\frac{2+\alpha(Q_N+iy)-n}{2}\right)\Gamma\left(\frac{\alpha \left(Q_N+iy\right)}{2}\right)}\right|\leq C \frac{(\sqrt{2}e^{\alpha })^{Q_N}}{Q_N^{Q_N+\alpha Q_N+\frac{1}{2}-\nu-\frac{n}{2}}}y^{Q_N-\frac{1}{2}}e^{\frac{\pi}{4}y+\frac{\pi\alpha }{2}y}.
\end{equation*}On the other hand, if  $v>0$, we have
\begin{equation*}
\left|\frac{1}{\sin (u+iv)}\right|=\left|\frac{2}{e^{iu-v}-e^{-iu+v}}\right|\leq \frac{2}{e^v - e^{-v}}.
\end{equation*}Therefore, if $ 1\ll Q_N\leq y$,
\begin{equation*}
\left|   f_{\alpha, \nu, n}(Q_N+iy; r)\right|\leq C \frac{(\sqrt{2}e^{\alpha})^{ Q_N}}{Q_N^{Q_N+\alpha Q_N+\frac{1}{2}-\nu-\frac{n}{2}}}y^{Q_N-\frac{1}{2}}e^{-\frac{3\pi}{4}y}\left(\frac{r}{2}\right)^{\alpha Q_N},
\end{equation*}and it follows that
\begin{equation*}\begin{split}
\left|\int_{|y|\geq Q_N}  f_{\alpha, \nu, n}(Q_N+iy; r) dy \right|\leq & C \frac{(\sqrt{2}e^{\alpha})^{ Q_N}}{Q_N^{Q_N+\alpha Q_N+\frac{1}{2}-\nu-\frac{n}{2}}}\Gamma\left(Q_N+\frac{1}{2}\right)\left(\frac{4}{3\pi}\right)^{Q_N+\frac{1}{2}}\left(\frac{r}{2}\right)^{\alpha Q_N}\\
&\leq C\frac{(\sqrt{2}e^{\alpha-1})^{ Q_N}}{Q_N^{\alpha Q_N+\frac{1}{2}-\nu-\frac{n}{2}}} \left(\frac{r}{2}\right)^{\alpha Q_N},
\end{split}\end{equation*}which tends to zero when $N\rightarrow \infty$, uniformly for $r$ in any compact subset of $[0, \infty)$. For $y\leq Q_N$,   \eqref{eq11_07_12_1} implies that\begin{equation*}
\left|\frac{1}{\sin\pi(Q_N+iy-\nu)}\right|\leq \left|\frac{1}{\sin \pi (Q_N -\nu )}\right| \leq \frac{1}{\sin \frac{\pi}{4}},\end{equation*}and
\begin{equation*}
\left|\frac{1}{\sin\frac{\pi\left( \alpha(Q_N+iy)-n\right)}{2}}\right|\leq \left|\frac{1}{\sin \frac{\pi\left( \alpha Q_N-n\right)}{2}}\right| \leq \frac{1}{\sin \frac{\pi \alpha}{8}}.\end{equation*}
Consequently,
\begin{equation*}
\left|\int_{|y|\leq Q_N}  f_{\alpha, \nu, n}(Q_N+iy; r) dy \right|\leq  C \frac{\left(\sqrt{2}(2e/\alpha)^{\alpha}e^{\left(\alpha+1\right)\frac{\pi}{4}}\right)^{Q_N}}{Q_N^{\alpha Q_N-\nu-\frac{n}{2}}}\left(\frac{r}{2}\right)^{\alpha Q_N},
\end{equation*}which tends to zero when $N\rightarrow \infty$, uniformly for $r$ in any compact subset of $[0, \infty)$.
\end{proof}

Now we want to investigate the analyticity of the series representation of $q_{\alpha,\nu, n}(r)$. First we have

\begin{proposition}\label{pro7}
Let $\alpha =a/b\in (0,2)$ and $\nu=e/f>0$ be rational numbers with $\text{gcd}\, (a,b)=\text{gcd}\,(e,f)=1$. If one of the conditions \\(i) $a$ is not divisible by $f$, \\ (ii) $n$ is an odd integer and $a/f$ is an even integer, \\(iii) $n$ is an even integer, $a$ and $f$ are both even integers and $a/f$ is an odd integer,\\ holds, then $q_{\alpha, \nu, n}(r)$ can be represented by:
\begin{equation}\label{eq11_06_10}
q_{\alpha,\nu, n}(r)=\frac{1}{2^n\pi^{\frac{n}{2}}\Gamma(\nu)}\left\{ \frac{r^{\alpha\nu-n}}{2^{\alpha\nu-n}}A_{1;\alpha,\nu,n}\left( r^{\alpha}\right)+A_{2;\alpha,\nu,n}\left(r^2\right)\right\},
\end{equation}
where $A_1(z)$ and $A_2(z)$ are entire functions given by
\begin{equation}\label{eq11_06_5}\begin{split}
A_{1;\alpha,\nu,n}(z)=&\sum_{l=0}^{\infty} \frac{(-1)^{l}}{l!}\Gamma(\nu+l)\frac{\Gamma\left(\frac{n-\alpha(\nu+l)}{2}\right)}{\Gamma\left(\frac{\alpha(\nu+l)}{2}\right)} \frac{z^{l}}{2 ^{ \alpha l }}\\A_{2; \alpha, \nu,n}(z)=&\frac{2}{\alpha}\sum_{j=0}^{\infty}\frac{(-1)^j}{j!}
\frac{\Gamma\left(\frac{n+2j}{\alpha}\right)}{\Gamma\left(\frac{n+2j}{2}\right)}\Gamma\left(\frac{\alpha \nu -n-2j}{\alpha}\right) \frac{z^j}{2^{2j}} .\end{split}
\end{equation}
\end{proposition}
\begin{proof}If the conditions (i) or (ii) holds, Proposition \ref{pro15} implies that $(\alpha,\nu)\notin\Lambda_n$.
Using  Propositions \ref{pro4} and \ref{pro6}, we only need to show that the right hand sides of the equations \eqref{eq11_06_5} define entire functions. Let
\begin{equation*}\begin{split}
A_{1, l} =& \frac{(-1)^{l}}{l!}\Gamma(\nu+l)\frac{\Gamma\left(\frac{n-\alpha(\nu+l)}{2}\right)}{\Gamma\left(\frac{\alpha(\nu+l)}{2}\right)}= \frac{(-1)^{l-1}}{l!}\frac{\Gamma(\nu+l)}{\Gamma\left(\frac{2+\alpha(\nu+l)-n}{2}\right)\Gamma\left(\frac{\alpha(\nu+l)}{2}\right)}\frac{\pi}{\sin \frac{\pi\left(\alpha(\nu+l)-n\right)}{2}},\\
A_{2,j}= &\frac{(-1)^j}{j!}
\frac{\Gamma\left(\frac{n+2j}{\alpha}\right)}{\Gamma\left(\frac{n+2j}{2}\right)}\Gamma\left(\frac{\alpha \nu -n-2j}{\alpha}\right)=
\frac{(-1)^{j-1}}{j!}
\frac{\Gamma\left(\frac{n+2j}{\alpha}\right)}{\Gamma\left(\frac{n+2j}{2}\right)}\frac{\pi }{\Gamma\left(\frac{\alpha-\alpha \nu +n+2j}{\alpha}\right)\sin \frac{\pi\left(n+2j-\alpha\nu\right)}{\alpha}}.\end{split}
\end{equation*}Given $l\in \hat{\mathbb{N}}$, there exists a unique $j_l$ such that
\begin{equation*}
\left|\frac{\alpha(\nu+l)-n}{2}-j_l\right|\leq \frac{1}{2}.
\end{equation*}Since $\frac{\alpha(\nu+l)-n}{2}-j_l \neq 0$,
\begin{equation*}
\left|\frac{\alpha(\nu+l)-n}{2}-j_l\right| =\left|\frac{a(e+fl)-bf(n+2j_l)}{2bf} \right| \geq \frac{1}{2bf}.
\end{equation*}Similarly, given $j\in \hat{\mathbb{N}}$, there exists a unique $l_j$ such that
\begin{equation*}
\left|\frac{n+2j-\alpha\nu}{\alpha}-l_j\right|\leq \frac{1}{2}.
\end{equation*}$\frac{n+2j-\alpha\nu}{\alpha}-l_j\neq 0$ implies that
\begin{equation*}
\left|\frac{n+2j-\alpha\nu}{\alpha}-l_j\right|=\left| \frac{bf(n+2j) -ae-al_j}{af}\right|\geq \frac{1}{af}.
\end{equation*} Therefore,
\begin{equation*}
\left|\frac{1}{\sin \frac{\pi\left(\alpha(\nu+l)-n\right)}{2}}\right|=\left|\frac{1}{\sin \pi\left(\frac{ \alpha(\nu+l)-n }{2}-j_l\right)}\right|\leq \frac{1}{\sin \frac{\pi}{2bf}}\end{equation*}and\begin{equation*}
\left|\frac{1}{\sin \frac{\pi \left(n+2j-\alpha\nu\right)}{\alpha}}\right|=\left|\frac{1}{\sin \pi \left(\frac{  n+2j-\alpha\nu }{\alpha}-l_j\right)}\right|\leq \frac{1}{\sin \frac{\pi}{af}}.
\end{equation*}
Together with
Stirling's formula \eqref{eq11_13_5}, we find that for $l$ and $j$ large enough,
\begin{equation*}
\frac{|A_{1, l+1}|}{|A_{1, l}|}\leq \frac{C_1}{l^{\alpha}\sin \frac{\pi}{2bf}}\xrightarrow{l\rightarrow \infty} 0
\hspace{1cm}\text{and}\hspace{1cm}
\frac{|A_{2, j+1}|}{|A_{2, j}|}\leq \frac{C_2}{j^2\sin \frac{\pi}{af}}\xrightarrow{j\rightarrow \infty} 0.
\end{equation*}Therefore, the two infinite series on the right hand sides of the equations \eqref{eq11_06_5} converge absolutely and uniformly on any compact subsets of $\mathbb{C}$.
\end{proof}

\begin{proposition}\label{pro8}
Let $\alpha =a/b\in (0,2)$ and $\nu=e/f>0$ be rational numbers with $a, b, e, f\in \mathbb{N}$ and $\text{gcd}\, (a,b)=\text{gcd}\,(e,f)=1$. If one of the conditions \\(i) $n$ is an odd integer and $a/f$ is an odd integer,\\(ii) $n$ is an even integer and $a/f$ is an even integer, \\(iii) $n$ is an even integer, $a/f$ is an odd integer and $a$ and $f$ are both odd integers,\\ holds, let $(j_0, l_0)$ be the smallest nonnegative integers satisfying $n+2j= \alpha(\nu+l), (j, l)\in \hat{\mathbb{N}}\times\hat{\mathbb{N}} $. Then $q_{\alpha, \nu, n}(r)$ can be represented by:
\begin{equation}\label{eq11_06_9}\begin{split}
q_{\alpha,\nu, n}(r)=&\frac{1}{2^n\pi^{\frac{n}{2}}\Gamma(\nu)}\Biggl\{ \frac{r^{\alpha\nu-n}}{2^{\alpha\nu-n}}A_{1;\alpha,\nu,n}\left( r^{\alpha}\right)+A_{2;\alpha,\nu,n}\left(r^2\right)\\&+ \frac{r^{2j_0}}{2^{2j_0}}A_{3;\alpha,\nu,n}(r^a)+\frac{r^{2j_0}}{2^{2j_0}}A_{4;\alpha,\nu,n}(r^{a})\log \frac{r}{2}\Biggr\},\end{split}
\end{equation}
where if $a$ is odd, $A_1(z)$, $A_2(z)$, $A_3(z)$ and $A_4(z)$ are entire functions given by
\begin{equation}\label{eq11_06_8}\begin{split}A_{1;\alpha,\nu,n}(z)=&\sum_{\substack{l\in \hat{\mathbb{N}};\; l \neq l_0+2bq, q\in\hat{\mathbb{N}}} } \frac{(-1)^{l}}{l!}\Gamma(\nu+l)\frac{\Gamma\left(\frac{n-\alpha(\nu+l)}{2}\right)}{\Gamma\left(\frac{\alpha(\nu+l)}{2}\right)}  \frac{z^l}{2^{ \alpha l }},\\A_{2;\alpha,\nu,n}(z)=&\frac{2}{\alpha}\sum_{\substack{j\in \hat{\mathbb{N}};\;j\neq j_0+ aq, q\in \hat{\mathbb{N}}}}\frac{(-1)^j}{j!}
\frac{\Gamma\left(\frac{n+2j}{\alpha}\right)}{\Gamma\left(\frac{n+2j}{2}\right)}\Gamma\left(\frac{\alpha \nu -n-2j}{\alpha}\right) \frac{z^j}{2^{2j}}, \\ A_{3;\alpha,\nu,n}(z)=&-\frac{2(-1)^{l_0+j_0}}{\alpha}\sum_{q=0}^{\infty}\frac{(-1)^{q}}
{(l_0+2qb)!\left(j_0+ aq \right)!}\frac{\Gamma(\nu+l_0+2bq)}{\Gamma\left(\frac{n}{2}+j_0+ aq \right)}\Biggl\{\psi(\nu+l_0+2bq)\\&-\psi(l_0+2bq+1)-\frac{a}{2b}
\psi\left(\frac{n}{2}+j_0+ aq \right) -\frac{a}{2b}\psi\left(j_0+ aq +1\right) \Biggr\} \frac{z^{2q}}{2^{2aq}} ,\\
 A_{4;\alpha,\nu,n}(z)=&- 2(-1)^{l_0+j_0} \sum_{q=0}^{\infty}\frac{(-1)^{q }}
{(l_0+2qb)!\left(j_0+ aq \right)!}\frac{\Gamma(\nu+l_0+2bq)}{\Gamma\left(\frac{n}{2}+j_0+aq\right)}  \frac{z^{2q}}{2^{2aq}}; \end{split}
\end{equation}whereas if $a$ is even,\begin{equation}\label{eq11_06_6}\begin{split}
A_{1;\alpha,\nu,n}(z)=&\sum_{\substack{l\in \hat{\mathbb{N}};\; l \neq l_0+bq, q\in\hat{\mathbb{N}}} } \frac{(-1)^{l}}{l!}\Gamma(\nu+l)\frac{\Gamma\left(\frac{n-\alpha(\nu+l)}{2}\right)}{\Gamma\left(\frac{\alpha(\nu+l)}{2}\right)}  \frac{z^l}{2^{ \alpha l }},\\A_{2;\alpha,\nu,n}(z)=&\frac{2}{\alpha}\sum_{\substack{j\in \hat{\mathbb{N}};\;j\neq j_0+\frac{aq}{2}, q\in \hat{\mathbb{N}}}}\frac{(-1)^j}{j!}
\frac{\Gamma\left(\frac{n+2j}{\alpha}\right)}{\Gamma\left(\frac{n+2j}{2}\right)}\Gamma\left(\frac{\alpha \nu -n-2j}{\alpha}\right) \frac{z^j}{2^{2j}}, \\ A_{3;\alpha,\nu,n}(z)=&-\frac{2(-1)^{l_0+j_0}}{\alpha}\sum_{q=0}^{\infty}\frac{(-1)^{q\left(b+\frac{a}{2}\right)}}
{(l_0+qb)!\left(j_0+\frac{aq}{2}\right)!}\frac{\Gamma(\nu+l_0+bq)}{\Gamma\left(\frac{n}{2}+j_0+\frac{aq}{2}\right)}\Biggl\{\psi(\nu+l_0+bq)\\&-\psi(l_0+bq+1)-\frac{a}{2b}
\psi\left(\frac{n}{2}+j_0+\frac{aq}{2}\right) -\frac{a}{2b}\psi\left(j_0+\frac{aq}{2}+1\right) \Biggr\} \frac{z^q}{2^{aq}},\\
 A_{4;\alpha,\nu,n}(z)=&- 2(-1)^{l_0+j_0} \sum_{q=0}^{\infty}\frac{(-1)^{q\left(b+\frac{a}{2}\right)}}
{(l_0+qb)!\left(j_0+\frac{aq}{2}\right)!}\frac{\Gamma(\nu+l_0+bq)}{\Gamma\left(\frac{n}{2}+j_0+\frac{aq}{2}\right)}  \frac{z^q}{2^{aq}}. \end{split}
\end{equation}

\end{proposition}\begin{proof}
The fact that $q_{\alpha, \nu, n}(r)$ can be represented in the form \eqref{eq11_06_9} follows from Propositions \ref{pro4}, \ref{pro15} and \ref{pro6}. The proof of the absolute convergence of the series $A_{1;\alpha,\nu,n}(z)$ and $A_{2;\alpha,\nu,n}(z)$, $A_{3;\alpha,\nu,n}(z)$ and $A_{4;\alpha,\nu,n}(z)$ follows the same as the proof of Proposition \ref{pro7} and the formula (see e.g. \cite{73}):
\begin{equation*}
\psi(x) = \log x -\frac{1}{2x} + O\left(\frac{1}{x^2}\right)\hspace{1cm}\text{as}\;\;x\rightarrow \infty.
\end{equation*}

\end{proof}

Propositions \ref{pro7} and \ref{pro8} show that when $\alpha$ and $\nu$ are both rational numbers,   $q_{\alpha,\nu, n}(r)$ can be represented in terms of  entire functions. In the case where $\alpha$ is a rational number and $\nu$ is an irrational number, we have
\begin{proposition}\label{pro9}
Let $\alpha=a/b\in (0, 2)$ be a rational number with $a, b\in \mathbb{N}$ and $(a,b)=1$. If $\nu>0$ is an irrational number,
 then $q_{\alpha,\nu, n}(r)$ can be written as the sum of two absolutely convergent series as in Proposition \ref{pro7}.

\end{proposition}
\begin{proof}
Notice that the conditions $\alpha$ is rational and $\nu$ is irrational imply that $(\alpha,\nu)\notin\Lambda_n$.    Given $\nu>0$ an irrational number, let $q_{\nu}$ be the unique positive integer such that
\begin{equation*}
\delta_{\nu}=\left|\nu -\frac{q_{\nu}}{a}\right|<\frac{1}{2a}.
\end{equation*} This implies that for all $q\in\mathbb{Z}$
\begin{equation}\label{eq11_10_4}
\left|\nu -\frac{q}{a}\right|\geq \min\left\{\delta_{\nu},\frac{1}{2a}\right\}:=\eta_{\nu}.
\end{equation} Given $l\in \hat{\mathbb{N}}$, there exists a unique $j_l$ such that
\begin{equation*}
\left|\frac{\alpha(\nu+l)-n}{2}-j_l\right|\leq \frac{1}{2}.
\end{equation*}\eqref{eq11_10_4} implies that
\begin{equation*}
\left|\frac{\alpha(\nu+l)-n}{2}-j_l\right|=\frac{\alpha}{2}\left| \nu - \frac{b(n+2j_l)-al}{a}\right|\geq \frac{\alpha\eta_{\nu}}{2}.
\end{equation*}Similarly, given $j\in \hat{\mathbb{N}}$, there exists a unique $l_j$ such that
\begin{equation*}
\left|\frac{n+2j-\alpha\nu}{\alpha}-l_j\right|\leq \frac{1}{2}.
\end{equation*}\eqref{eq11_10_4} implies that
\begin{equation*}
\left|\frac{n+2j-\alpha\nu}{\alpha}-l_j\right|=\left|\nu - \frac{b(n+2j)-al_j}{a}\right|\geq \eta_{\nu}.
\end{equation*}The rest of the proof is the same as Proposition \ref{pro7}.

\end{proof}

When $\nu$ is a rational number but $\alpha$ is irrational,  the situation is more complicated. We have
to introduce the concept of Liouville numbers as in \cite{10, 18, 20, 19, 17}. Recall that an irrational number $\beta$ is called a Liouville number if for all $m=2, 3, \ldots$, there exists a rational number $p/q, p, q\in \mathbb{Z}$ such that
\begin{equation*}
\left| \beta -\frac{p}{q}\right|\leq \frac{1}{q^m}.
\end{equation*}By Liouville theorem (see e.g. \cite{74}), all Liouville numbers are transcendental and the set of Liouville numbers, denoted by $L$, has Lebesgue measure zero.
\begin{proposition}\label{pro10}
Let $\nu =e/f>0$ be a rational number with $e, f\in \mathbb{N}, (e, f)=1$. If
 $\alpha$ is an irrational number but not a Liouville number, then $q_{\alpha,\nu, n}(r)$ can be written as the sum of two absolutely convergent series as in Proposition \ref{pro7}.
\end{proposition}
\begin{proof}
Notice that the conditions on $\alpha$ and $\nu$ imply that $(\alpha,\nu)\notin \Lambda_n$.
 Since $\alpha$ is  not a Liouville number, there exists an integer $m_{\alpha}$ such that for all $m\geq m_{\alpha}$, and for all $p\in \hat{\mathbb{N}}, q\in \mathbb{N}$,
\begin{equation}\label{eq11_10_3}\left| \alpha-\frac{p}{q}\right|\geq \frac{1}{q^m}.\end{equation}
Fix $l \in \hat{\mathbb{N}}$. There exists a unique $j_l\in \mathbb{Z}$, such that
\begin{equation*}
\left|\frac{\alpha(\nu+l)-n}{2}-j_l\right|\leq \frac{1}{2}.
\end{equation*} By \eqref{eq11_10_3}, we have\begin{equation*}
 \left|\frac{\alpha(\nu+l)-n}{2}-j_l\right|=\frac{(\nu+l)}{2}\left|\alpha -\frac{f(n+2j_l)}{e+fl}\right| \geq \frac{1}{2f^{m_{\alpha}}(\nu+l)^{m_{\alpha}-1}}.
\end{equation*} This implies that
\begin{equation}\label{eq11_10_1}
\left|\frac{1}{\sin \frac{\pi(\alpha(\nu+l)-n)}{2}}\right|=\left|\frac{1}{\sin \pi\left(\frac{ \alpha(\nu+l)-n}{2}-j_l\right)}\right|\leq \frac{2f^{m_{\alpha}}(\nu+l)^{m_{\alpha}-1}}{\pi}.
\end{equation}On the other hand, given $j\in \hat{\mathbb{N}}$, let $l_j$ be the unique integer such that
\begin{equation*}
-\frac{1}{\alpha}<-\frac{1}{2}\leq \frac{\left(n+2j \right)}{\alpha}-(\nu+l_j) \leq \frac{1}{2}< \frac{1}{\alpha}.
\end{equation*}This implies that
\begin{equation*}
\frac{1}{\nu + l_j}> \frac{\alpha}{n+2j+1}.
\end{equation*}Consequently,
\begin{equation*}\begin{split}
\left| \frac{\left(n+2j-\alpha\nu\right)}{\alpha}-l_j\right|&=\frac{\nu+l_j}{\alpha}\left|\frac{f(n+2j)}{e+fl_j}-\alpha\right|\\
&\geq \frac{1}{\alpha f^{m_{\alpha}}(\nu +l_j)^{m_{\alpha}-1}}\geq \frac{\alpha^{m_{\alpha}-2}}{f^{m_{\alpha}}(n+2j+1)^{m_{\alpha}-1}}, \end{split}
\end{equation*}and
\begin{equation}\label{eq11_10_2}
\left|\frac{1}{\sin \frac{\pi(n+2j-\alpha\nu)}{\alpha}}\right|\leq \frac{f^{m_{\alpha}}(n+2j+1)^{m_{\alpha}-1}}{\pi \alpha^{m_{\alpha}-2}}.
\end{equation}If   $m_{\alpha}$ is larger than $\alpha$, the ratio $|A_{1, l+1}|/|A_{1, l}|$ or $|A_{2, j+1}|/|A_{2, j}|$ may become unbounded as $l\rightarrow \infty$ or $j\rightarrow \infty$. Therefore, instead of considering the ratio, we consider the roots $ |A_{1, l}|^{\frac{1}{l}}$ and $ |A_{2, j}|^{\frac{1}{j}}$. Stirling's formula \eqref{eq11_13_5}, \eqref{eq11_10_1} and \eqref{eq11_10_2} imply that
\begin{equation*}
|A_{1, l}|^{\frac{1}{l}}\leq C_1 \frac{l^{m_{\alpha}/l}}{l^{\alpha}}\xrightarrow{l\rightarrow \infty} 0\hspace{1cm}\text{and}\hspace{1cm}|A_{2, j}|^{\frac{1}{j}}\leq C_2 \frac{j^{m_{\alpha}/j}}{j^2}\xrightarrow{j\rightarrow \infty} 0.
\end{equation*}This shows that both the series $A_{1;\alpha,\nu,n}(z)$ and $A_{2;\alpha,\nu,n}(z)$ are absolutely and uniformly convergent for all $z \in \mathbb{C}$.

\end{proof}
Recall that an entire function $E(z)=\sum_{k=0}^{\infty} E_k z^k$ is said to be of order $\rho$ if and if
\begin{equation*}
\rho=\limsup_{k\rightarrow \infty}\frac{k\log k}{\log 1/|E_k|}.
\end{equation*}It is easy to deduce from the proof of Proposition \ref{pro10} that when $\alpha,\nu$ are both rational numbers, or if one of $\alpha$ or $\nu$ is a rational number, and the other is a non-Liouville irrational number, then the entire functions $A_{1;\alpha,\nu,n}(z)$ and $A_{2;\alpha,\nu,n}(z)$ have order $1/\alpha$ and $1/2$ respectively. In the case where $\alpha=a/b$ and $\nu=e/f$ satisfy the conditions (i), (ii) or (iii) of Proposition \ref{pro8}, the entire functions $A_{3;\alpha,\nu,n}(z)$ and $A_{4;\alpha,\nu,n}(z)$ have order $1/a$.

Although Propositions \ref{pro7}, \ref{pro8} and \ref{pro9} are proved for $\alpha\in (0,2)$, one can show that by naively extending the results to $\alpha=2$, one obtains  exactly the results for $\alpha=2$ given in Proposition \ref{pro11}.

At this point, one may wonder whether there exist any values of $\alpha$ and $\nu$ where both the series $A_{1;\alpha,\nu,n}(z)$ and $A_{2;\alpha,\nu,n}(z)$ are divergent. As was proved in \cite{18}, when $\nu=1$, there is a dense subset of $\alpha\in (0,2)$ where both the series $A_{1;\alpha,\nu,n}(z)$ and $A_{2;\alpha,\nu,n}(z)$ are divergent. Now we extend the result to general $\nu>0$ which can also be regarded as the generalization of a result in \cite{19} for $n=1$.
\begin{proposition}Given a rational number $\nu = e/f>0$, there exists a dense subset of $\alpha\in (0,2)$ where both the series $A_{1;\alpha,\nu,n}(z)$ and $A_{2;\alpha,\nu,n}(z)$ are divergent. This implies that there is a dense subset of $(\alpha,\nu)\in (0,2)\times (0, \infty)$ where  both the series $A_{1;\alpha,\nu,n}(z)$ and $A_{2;\alpha,\nu,n}(z)$ are divergent.
\end{proposition}
\begin{proof}The proof follows closely the ideas of \cite{18, 19}. Given $\nu=e/f$, define a sequence of positive integers $\{\sigma_k\}_{k=1}^{\infty}$ by $\sigma_1=1+2f$ and $\sigma_{k+1}=(1+2f)^{2\sigma_{k}}$.
Let $\Delta$ be the set of sequences $\{\delta_m\}_{m=1}^{\infty}$ satisfying the conditions:

(i) $\delta_m$ is equal to 0 or 1.

(ii) $\delta_m=1$ implies that $m\in \{\sigma_k\}_{k=1}^{\infty}$.

(iii) $\delta_m=1$ for infinitely many $m$.\\
Define $\Omega$ to be the set of real numbers of the form $$\sum_{m=1}^{\infty}\frac{\delta_m}{(2f+1)^{m}},$$
and let
$\Lambda$   be the set of rational numbers of the form
\begin{equation*}
(1+2f)^s\sum_{m=0}^{t}\frac{a_m}{(1+2f)^{m}}, \hspace{1cm} s, t\in \hat{\mathbb{N}}, a_m \in \{0, 1, 2, \ldots, 2f\}.
\end{equation*}
Finally define $E$ to be the set
\begin{equation*}
E=\left\{ \alpha\in (0,2)\,:\, \alpha\nu =\frac{e\alpha}{f}=x + y, \;x\in\Lambda, y\in \Omega\right\}.
\end{equation*}Clearly $E$ is a dense subset of $(0,2)$. Given $\alpha\in E$, we claim that for any $r>0$, there exists a subsequence of positive integers $l_k$ and a subsequence of integers $j_k$ such that
$$|A_{1, l_k}z^{l_k}|\xrightarrow{k\rightarrow \infty} \infty\hspace{1cm}\text{and}\hspace{1cm}|A_{2, j_k}z^{j_k}|\xrightarrow{k\rightarrow \infty} \infty$$for any $z\neq 0$. This will imply that    the series $A_{1;\alpha,\nu,n}(z)$ and $A_{2;\alpha,\nu,n}(z)$ are divergent.
Given $\alpha\in E$, $\alpha\nu$ can be written uniquely as
\begin{equation*}
\alpha\nu = h+  \sum_{m=1}^{t} \frac{a_m}{(2f+1)^m}+\sum_{m=t+1}^{\infty}\frac{\delta_m}{(2f+1)^{m}}, \hspace{1cm}h\in \hat{\mathbb{N}},\; a_m \in \{0, 1,\ldots, 2f\}.
\end{equation*}
Define the sequence $\{\eta_{k}\}_{k=1}^{\infty}$ so that $\{m>t\,:\,\delta_m=1\}=\{\eta_{k}\}_{k=1}^{\infty}$. Notice that $\eta_1\geq t+1$ and $$\eta_{k+1}\geq (2f+1)^{2\eta_k}\geq 9\eta_k.$$
For any $k\in \mathbb{N}$, we have
\begin{equation*}
0< \alpha\nu -  \left\{h+  \sum_{m=1}^{t} \frac{a_m}{(2f+1)^m}+\sum_{m=1}^{k}\frac{1}{(2f+1)^{\eta_m}}\right\}=\sum_{m=k+1}^{\infty}\frac{1}{(2f+1)^{\eta_m}}<\frac{1}{(2f+1)^{\eta_{k+1}-1}}.
\end{equation*}This implies that
\begin{equation*}\begin{split}
&\alpha\nu (2f+1)^{\eta_k}-\left( h(2f+1)^{\eta_k}+\sum_{m=1}^t a_m (2f+1)^{\eta_k-m}+\sum_{m=1}^k (2f+1)^{\eta_k-\eta_m}\right)\\& < \frac{1}{(2f+1)^{\eta_{k+1}-\eta_k-1}} <\frac{1}{(2f+1)^{\eta_{k+1}/2}}.
\end{split}\end{equation*}Notice that since $(2f+1)$ is odd, the integer
\begin{equation*}
p_k = h(2f+1)^{\eta_k}+\sum_{m=1}^t a_m (2f+1)^{\eta_k-m}+\sum_{m=1}^k (2f+1)^{\eta_k-\eta_m}
\end{equation*}is odd if and only if
\begin{equation*}
h+\sum_{m=1}^t a_m + k \hspace{1cm}\text{is odd}.
\end{equation*} If $h+\sum_{m=1}^t a_m-n $ is even, $p_{2k}-n$ is even.   In this case, we define $l_k $ and $j_k$ by $$l_k=\frac{e\left[(2f+1)^{\eta_{2k}}-1\right]}{f},\hspace{1cm}j_k =\frac{p_{2k}-n}{2}.$$If $h+\sum_{m=1}^t a_m-n $ is odd, $p_{2k-1}-n$ is even, and we define $l_k $ and $j_k$ by $$l_k=\frac{e\left[(2f+1)^{\eta_{2k-1}}-1\right]}{f},\hspace{1cm}j_k =\frac{p_{2k-1}-n}{2}.$$In the first case where $h+\sum_{m=1}^t a_m-n $ is even, we have
\begin{equation*}\begin{split}
\left|\frac{\alpha(\nu+l_k)-n}{2}-j_k\right| =& \frac{\alpha\nu (2f+1)^{\eta_{2k}}-n-p_{2k}+n}{2}=\frac{\alpha\nu (2f+1)^{\eta_{2k}}- p_{2k} }{2}\\<&\frac{1}{2(2f+1)^{\eta_{2k+1}/2}}
\end{split}\end{equation*}and\begin{equation*}
\begin{split}
\left|\frac{n+2j_k-\alpha\nu}{\alpha} -l_k \right|=\frac{2}{\alpha}\left|\frac{\alpha(\nu+l_k)-n}{2}-j_k\right|<\frac{1}{\alpha(2f+1)^{\eta_{2k+1}/2}}.
\end{split}
\end{equation*}
These imply that
\begin{equation*}
\left|\frac{1}{\sin \frac{\pi(\alpha(\nu+l_k)-n)}{2}}\right|=\left|\frac{1}{\sin \pi\left(\frac{\alpha(\nu+l_k)-n}{2}-j_k\right)}\right|\geq \frac{2(2f+1)^{\eta_{2k+1}/2}}{\pi}
\end{equation*}and
\begin{equation*}
\left|\frac{1}{\sin \frac{\pi(n+2j_k-\alpha\nu)}{\alpha}}\right|=\left|\frac{1}{\sin \pi\left(\frac{\alpha(\nu+l_k)-n-2j_k}{\alpha} \right)}\right|\geq  \frac{\alpha (2f+1)^{\eta_{2k+1}/2}}{\pi}.
\end{equation*}
Together with Stirling's formula \eqref{eq11_13_5}, we find that for $k\gg 1$,
\begin{equation*}
|A_{1, l_k} z^{l_k}|\geq \frac{C(2f+1)^{\eta_{k+1}/2}}{l_k^{3l_k}}|z|^{l_k}\;\;\text{and}\;\;|A_{2, j_k}z^{j_k}|\geq \frac{C (2f+1)^{\eta_{2k+1}/2}}{j_k^{3j_k}}|z|^{j_k}.
\end{equation*}Notice that $l_{k}\leq \nu(2f+1)^{\eta_{2k}}$ and $j_k \leq (h+t+k)(2f+1)^{\eta_{2k}}$. This implies that for $k$ large enough
\begin{equation*}\begin{split}
\frac{(2f+1)^{\eta_{2k+1}/2}}{l_k^{3l_k}}|z|^{l_k}\geq &\frac{(2f+1)^{\eta_{2k+1}/2}}{\left(\frac{\nu^3}{\sqrt{z}}(2f+1)^{3\eta_{2k}}\right)^{\nu(2f+1)^{\eta_{2k}}}}\\
\geq &\frac{(2f+1)^{\eta_{2k+1}/2}}{\left( (2f+1)^{4\eta_{2k}}\right)^{\nu(2f+1)^{\eta_{2k}}}}=(2f+1)^{\frac{\eta_{2k+1}}{2}-4\nu\eta_{2k}(2f+1)^{\eta_{2k}}}.
\end{split}\end{equation*}Now   \begin{equation*}\begin{split}
\frac{\eta_{2k+1}}{2}-4\nu\eta_{2k}(2f+1)^{\eta_{2k}}\geq &\frac{1}{2}(2f+1)^{2\eta_{2k}}-4\nu\eta_{2k}(2f+1)^{\eta_{2k}}\\
=& (2f+1)^{\eta_{2k}}\left( \frac{(2f+1)^{\eta_{2k}}}{2}- 4\nu \eta_{2k}\right)\xrightarrow{k\rightarrow \infty} \infty.\end{split}
\end{equation*}This shows that $|A_{1, l_k}z^{l_k}|\rightarrow \infty$ as $k\rightarrow \infty$. The proof for $|A_{2, j_k}z^{j_k}|\rightarrow \infty$ as $k\rightarrow \infty$ is similar. For the case where $h+\sum_{m=1}^t a_m-n $ is odd, the proof is the same. This concludes the assertion of the proposition.
\end{proof}

\end{document}